\documentclass[12pt]{amsart}

\calclayout

\usepackage{amssymb}
\usepackage[matrix,arrow,curve]{xy}
\usepackage{epsfig}

\usepackage{pstricks}
\usepackage{epic}
\usepackage{eepic}
\newtheorem{thm}{Theorem}[section]
\newtheorem{Theorem}[thm]{Theorem}

\newtheorem*{thma}{Theorem A}
\newtheorem*{thmb}{Theorem B}
\newtheorem*{thmc}{Theorem C}
\newtheorem{Lemma}[thm]{Lemma}

\newtheorem{Proposition}[thm]{Proposition}
\newtheorem{Corollary}[thm]{Corollary}
\theoremstyle{definition}
\newtheorem{Example}[thm]{Example}

\newtheorem{Remark}[thm]{Remark}

\newtheorem{ccote}[thm]{}

\newcommand{\preu}{\begin{proof}}
\newcommand{\preusp}[1]{\begin{proof}[Proof #1]}
\newcommand{\ppreu}[1]{\begin{proof}[Proof of #1]}

\newcommand{\mancqfd}{\unskip\kern 6pt\penalty 500
\raise -0pt\hbox{\vrule\vbox to7.3pt{\hrule width
6.66pt\vfill\hrule}\vrule}\smallskip}

\newcommand{\cqfd}{\end{proof}}
\newcommand{\proref}[1]{Proposition~\ref{#1}}
\newcommand{\remref}[1]{Remark~\ref{#1}}
\newcommand{\lemref}[1]{Lemma~\ref{#1}}
\newcommand{\corref}[1]{Corollary~\ref{#1}}
\newcommand{\thref}[1]{Theorem~\ref{#1}}
\newcommand{\exref}[1]{Example~\ref{#1}}
\newcommand{\secref}[1]{Section~\ref{#1}}

\newcommand{\bbd}{{\mathbb{D}}}
\newcommand{\bbs}{{\mathbb{S}}}

\newcommand{\bbr}{{\mathbb{R}}}
\newcommand{\bbt}{{\mathbb{T}}}
\newcommand{\bbc}{{\mathbb{C}}}

\newcommand{\bbz}{{\mathbb{Z}}}

\newcommand{\cala}{{\mathcal A}}

\newcommand{\calc}{{\mathcal C}}
\newcommand{\cald}{{\mathcal D}}
\newcommand{\cale}{{\mathcal E}}
\newcommand{\calf}{{\mathcal F}}

\newcommand{\cali}{{\mathcal I}}

\newcommand{\calo}{{\mathcal O}}

\newcommand{\calt}{{\mathcal T}}
\newcommand{\calu}{{\mathcal U}}
\newcommand{\calv}{{\mathcal V}}
\newcommand{\calw}{{\mathcal W}}
\newcommand{\calz}{{\mathcal Z}}

\newcommand{\bfz}{{\bf Z}}

\newcommand{\pcirc}{\kern .7pt {\scriptstyle \circ} \kern 1pt}

\newcommand{\hfl}[1]{\buildrel{#1}\over{\longrightarrow}}
\renewcommand{\eqref}[1]{(\ref{#1})}
\newcounter{exo}

\newcommand{\mun}{{-1}}
\newcommand{\sk}[1]{\vskip #1 mm}

\newcommand{\llangle}[2]{\langle #1 ,#2 \rangle}

\newcommand{\onto}{\to\kern-7.5pt\to}
\newcommand{\donto}{\downarrow\kern -7.92pt\raisebox{-0.6ex}{$\downarrow$}}

\newcommand{\algt}{\mathfrak{t}}

\newcommand{\algl}{\mathfrak{l}}

\newcommand{\scr}{\scriptscriptstyle}

\newcommand{\intb}{\,\raisebox{7.5pt}{$\scr\circ$}\kern -5pt}

\newcommand{\bung}{{\rm Bun}^G}
\newcommand{\nbungt}{{\rm Bun}^G_\Gamma}

\newcommand{\sbuns}{{\rm SBun}^{S^1}}
\newcommand{\bungs}{{\rm Bun}^{S^1}_\Gamma}
\newcommand{\buns}{{\rm Bun}^{S^1}}
\newcommand{\nsbungt}{{\rm SBun}^G_\Gamma}
\newcommand{\rep}{{\rm Rep}^G}

\newcommand{\reps}{{\rm Rep}^{S^1}}
\newcommand{\areps}{{\rm Rep}^{S^1}_{\scriptstyle\rm cell}}

\newcommand{\tilrep}{{\widetilde{\rm  Rep}}{}^G}

\newcommand{\arep}{{\rm Rep}^G_{\scriptstyle\rm cell}}

\newcommand{\carep}{\overline{\rm Rep}^{\,G}_{\scriptstyle\rm cell}}

\newcommand{\spl}{split }
\newcommand{\gagrou}{$(\Gamma,A)$-groupoid}
\newcommand{\cgagrou}{cellular $(\Gamma,A)$-groupoid}
\newcommand{\pgagrou}{proper $(\Gamma,A)$-groupoid}
\newcommand{\ctgrou}{cellular $(\bbt,A)$-groupoid}
\renewcommand{\hom}{{\rm Hom\,}}
\newcommand{\chom}{\overline{{\rm Hom}}\,}
\renewcommand{\varpi}{{\scriptstyle\Pi}}
\newcommand{\ciat}{\cali}
\newcommand{\repcat}{{\rm Rep}^G(\ciat)}

\newcommand{\repi}{{\rm Rep}^G(\cali)}

\newcommand{\isor}{{\bf\Phi}}

\DeclareMathOperator{\krep}{KRep}
\DeclareMathOperator{\sbun}{SBun}
\DeclareMathOperator{\repr}{Rep}

\newcommand{\lcom}{locally compact}
\newcommand{\loma}{weakly locally maximal}
\newcommand{\rloma}{locally maximal}
\newcommand{\sloma}{locally maximal}
\newcommand{\Sloma}{Locally maximal}

\newcommand{\lie}{{\rm Lie}}

\renewcommand{\:}{\colon}
\newcommand{\vv}{\,|\,}
\renewcommand{\labelenumi}{(\roman{enumi})} 
\newcommand{\grille}{%
\psline[linecolor=black](-4,0)(5,0)
\psline[linecolor=black](0,-4)(0,3)
\multiput(-4,0)(1,0){10}{\pscircle*[linecolor=black]{0.1}}%
\multiput(-4,1)(1,0){10}{\pscircle*[linecolor=black]{0.1}}%
\multiput(-4,2)(1,0){10}{\pscircle*[linecolor=black]{0.1}}%
\multiput(-4,3)(1,0){10}{\pscircle*[linecolor=black]{0.1}}%
\multiput(-4,-0)(1,0){10}{\pscircle*[linecolor=black]{0.1}}%
\multiput(-4,-1)(1,0){10}{\pscircle*[linecolor=black]{0.1}}%
\multiput(-4,-2)(1,0){10}{\pscircle*[linecolor=black]{0.1}}%
\multiput(-4,-3)(1,0){10}{\pscircle*[linecolor=black]{0.1}}%
\multiput(-4,-4)(1,0){10}{\pscircle*[linecolor=black]{0.1}}%
}

\begin{document}
\title{Equivariant bundles and isotropy representations
}
\author{Ian HAMBLETON}
\address{Department of Mathematics \& Statistics
\newline\indent
McMaster University
\newline\indent
Hamilton, ON L8S 4K1, Canada}
\email{ian@math.mcmaster.ca}
\author{Jean-Claude HAUSMANN}
\address{Section de Math\'ematiques
\newline\indent
Universit\'e de Gen\`eve,
 B.P. 240
\newline\indent
CH-1211 Gen\`eve 24,
Switzerland}
\email{Jean-Claude.Hausmann@unige.ch}
\date{May 28, 2009}
\thanks{\hskip -11pt  Research partially supported by NSERC Discovery Grant A4000. The authors would like to thank the Max Planck Institut f\"ur Mathematik for its hospitality and support while working on this paper.}

\begin{abstract}
We  introduce a new  
construction, the \emph{isotropy groupoid}, to organize the orbit  
data for split $\Gamma$-spaces.  We show that equivariant principal $G$-bundles over split $\Gamma$-CW complexes $X$ can be effectively classified by means of  
representations of their isotropy groupoids.  For instance, 
if the quotient complex  
$A=\Gamma\backslash X$ is a graph, with all edge stabilizers toral  
subgroups of $\Gamma$, we obtain a purely combinatorial  
classification of bundles with structural group $G$ a compact  
connected Lie group. 
If $G$ is abelian, our approach gives combinatorial and geometric 
descriptions of some
results of Lashof-May-Segal \cite{LMS}
and Goresky-Kottwitz-MacPherson \cite{GKM}.
\end{abstract}
\maketitle

\section*{Introduction}\label{intro}
In this paper we continue our study of equivariant principal bundles 
via isotropy representations (see \cite{hhausmann2}). 
If $\Gamma$ and $G$ are  topological groups, then 
a $\Gamma$-equivariant principal $G$-bundle is a locally trivial,
principal $G$-bundle $p\colon E\to X$ such that $E$ and $X$ 
are left $\Gamma$-spaces. The projection map $p$ is $\Gamma$-equivariant 
and $\gamma(e\cdot g)= (\gamma e)\cdot g$, where $\gamma \in \Gamma$ and
$g\in G$ acts on $e\in E$  by the principal action.
Equivariant principal bundles, and their natural generalizations, were 
studied by T.~E.~Stewart \cite{S},
T.~tom Dieck \cite{D1}, \cite[I\,(8.7)]{D2}, R.~Lashof
\cite{L2}, \cite{L3} together with P.~May \cite{LM} and G.~Segal \cite{LMS}.

The \emph{isotropy representation} at a  point $x \in X$ 
is the homomorphism $\alpha_{x}\colon\Gamma_x\to G$
defined by the formula 
$$\gamma\cdot \tilde x = \tilde x \cdot \alpha_x(\gamma)$$
where $\tilde x \in p^{-1}(x)$. The homomorphism $\alpha_x$ is independent
of the choice of $\tilde x$ up to conjugation in $G$. Here $\Gamma_x$ denotes the isotropy subgroup or stabilizer of $x\in X$.

The use of isotropy representations is particularly effective when 
the projection 
$\pi\: X \to \Gamma\backslash X\xrightarrow{\approx} A$ has a section 
$\varphi\: A \to X$.
We call the triple 
$(X,\pi,\varphi)$ a \emph{split} $\Gamma$-space over $A$. A natural source of 
examples is symplectic toric manifolds (see (\ref{ex-toric})), 
where $A$ is the moment polytope.
Under reasonable assumptions, a split $\Gamma$-space over $A$ 
is uniquely determined by its \emph{isotropy groupoid} 
$$
\cali:=\{(\gamma,a) \in \Gamma\times A\vv  \gamma\in\Gamma_{\varphi(a)}\}
$$
(see Proposition \ref{recon-pro}). 
 A $\Gamma$-equivariant principal $G$-bundle $\eta:=(E\xrightarrow{p} X)$  is called \emph{split} if
the pull-back $\varphi^*(\eta)$ is a trivial bundle. The isotropy representations of 
$\eta$ then produce a continuous groupoid representation
of $\cali$ in $G$ which is well defined up to conjugation by ${\rm Map}(A,G)$. We denote by 
$$
\repi=\hom(\cali,G)\big/{\rm Map}(A,G)
$$ 
the space of conjugacy classes of such groupoid representations. 
The first part of this paper is devoted to proving the following general classification theorem.

\begin{thma} Suppose that $\Gamma$ and $G$ are compact Lie groups. Let $X$ be split $\Gamma$-space over $A$ with isotropy groupoid $\cali$. Assume that $A$ is locally compact, and that $\cali$
is locally maximal. 
Then the equivalence classes of split $\Gamma$-equivariant $G$-bundles over $X$ are in bijection with $\repi$.
\end{thma}

The relevant definitions are given in Section \ref{nspeqbd}: see \ref{nspeqbd-defi} for equivariant bundles, \ref{S-isor} for the notion of a locally maximal isotropy groupoid,
and the proof of Theorem~A is given in \ref{S-classi}. In our applications
we will assume that $X$ is a $\Gamma$-CW-complex, equipped with a splitting 
$\varphi\: A \to X$ such that $\Gamma_{\varphi(a)}$ is constant on each open cell 
of its quotient CW-complex $A$. This property doesn't always hold (see 
\remref{gammaCW}), but it seems a natural assumption. We call the resulting isotropy groupoids \emph{cellular} (see \ref{stru}). A cellular groupoid is a combinatorial object and, when $\Gamma$
is discrete, it is 
a particular case of a \emph{developable simple complex of groups}  as considered by M.~Bridson and A.~Haefliger \cite{BH}. Cellular groupoids whose stalks $\cali_a$ are compact 
Lie groups are called
\emph{proper} groupoids. They arise from proper actions of 
a Lie group $\Gamma$, as studied for example in \cite{Lu-O} for $\Gamma$ discrete. 
In Theorem \ref{Th-classi-CW} we extend Theorem A to the classification
of split $\Gamma$-equivariant bundles over a proper groupoid.

In a second part of this
paper, we describe some approaches to computing $\repi$ assuming $\cali$ is cellular.
There is a corresponding notion of
\emph{cellular representations}, meaning those which are constant on the open cells of $A$,
whose classes modulo conjugation by a fixed element of $G$ form a set
denoted by $\arep(\cali)$ (see \ref{algrep}). The cellular representations are also purely combinatorial, and for $A$ a regular $CW$-complex, $\arep(\cali)$ is determined by restriction to the $0$-skeleton
of $A$ (see \proref{res01skel}).
We consider Theorem A to be an effective method of classifying equivariant bundles whenever $\repi$ can be reduced to $\arep(\cali)$. We therefore study the
natural map $\arep(\cali)\to\repi$, which turns out to be injective (\proref{jinjective}), but not surjective in general (see \eqref{kappanotonto}). It is however bijective
when $G$ is compact abelian (\proref{ijkabelian}), or when $A$ is a tree
(\proref{jinjective-cor}).

We next consider the case where $A$ is a graph. Here it is useful to define a related object $\carep(\cali)$, 
by allowing conjugation of cellular homomorphisms over each cell of $A$ independently.
It turns out that there exists a natural map $\imath\colon \repi\to\carep(\cali)$. 
In \thref{NTgraph}, we study this map for $G$ a compact connected Lie group. 
A sample application of \thref{NTgraph} is given by the following:

\begin{thmb}
Let $\Gamma$ and $G$ be a compact Lie groups, with $G$ connected, and suppose that $A$ is a graph. If $\cali$ is a cellular groupoid with all edge stabilizers torus subgroups of $\Gamma$, then the map  $\imath\colon \repi\to\carep(\cali)$ is a bijection.
\end{thmb}

Recall that there exists a classifying space $B(\Gamma, G)$  for $\Gamma$-equivariant principal $G$-bundle \cite{D1}, 
so the classification of equivariant bundles in
particular cases can also be approached by studying the $\Gamma$-equivariant
homotopy type of $B(\Gamma, G)$. 
If the structural group $G$ of the bundle
is {\it abelian}, then the main result of Lashof, May and Segal \cite{LMS} states that equivariant
bundles over a $\Gamma$-space $X$ are classified by the ordinary
homotopy classes of maps $[X \times_{\Gamma} E\Gamma, BG]$. For non-abelian structural groups, it appears that the natural map
$\theta\colon [X, B(\Gamma, G)]_\Gamma \to [X \times_{\Gamma} E\Gamma, BG]$ 
misses a lot of information, and our results could be interpreted as studying $\theta^{-1}(\bullet)$. 

Our isotropy groupoid $\cali$ has a classifying space 
$B\cali$ constructed by Haefliger \cite[p.~140]{Hae1}. We observe that
$B\cali \simeq X \times_{\Gamma} E\Gamma$ when $\cali$ is cellular. 
In our language, the result of \cite{LMS} implies that the natural map
$B\colon \repi \to [B\cali, BG]$ is injective for $G$ compact abelian. More generally, we show in \corref{mapB}:
\begin{thmc} Let $G$ and $\Gamma$ be compact Lie groups, with $G$ abelian. 
Let $X$ be a \spl $\Gamma$-space over $A$
with cellular isotropy groupoid $\ciat$. Suppose that $H^1(A;\pi_0(G))=H^2(A;\bbz)=0$. Then the map
$B\:\repi\to [B\cali,BG]$ is a bijection.
\end{thmc}
In \ref{ktheory} we point out the connection between our classification results and equivariant $K$-theory. Finally, in \ref{eqcoh}, we compare our results
with some classical theorems in equivariant cohomology, due to 
Chang-Skjelbred \cite{CS} and  Goresky-Kottwitz-MacPherson  \cite{GKM}.

\tableofcontents

\section{Preliminaries}\label{preli}

Most of this section contains folklore facts about compact Lie groups.
Let $K$ and $G$ be topological groups. We denote by $\hom(K,G)$ the space of continuous homomorphisms from $K$ to $G$, endowed with the compact-open topology.
Two homomorphisms $\alpha_1,\alpha_2\in\hom(K,G)$ are called {\it conjugate} if there
exists $g\in G$ such that $\alpha_2(\gamma)=g^\mun\alpha_1(\gamma)g$ for
all $\gamma\in K$. We denote by $\chom(K,G)$ the space of conjugacy classes,
endowed with the quotient topology.

\begin{Lemma}\label{condiscret-bis}
Let $K$ and $G$ be compact Lie groups. Then the space $\chom(K,G)$
is totally disconnected.
\end{Lemma}

\preu A bi-invariant Riemannian metric on $G$ gives rise to a
bi-invariant distance $d$ on $G$. The uniform convergence distance
on $\hom(K,G)$:
$$
d(\alpha,\beta)=\max_{k\in K}d(\alpha(k),\beta(k)) \ .
$$
induces the compact-open topology.
Let us denote by $\bar\alpha,\bar\beta\in\chom(K,G)$ the conjugacy classes
of $\alpha$ and $\beta$. One checks that the formula
$$
\bar d(\bar\alpha,\bar\beta)=\min_{g\in G}d(g\alpha g^\mun,\beta)=
\min_{g,h\in G}d(g\alpha g^\mun,h\beta h^\mun)
$$
defines a distance on $\chom(K,G)$. Since $\bar d(\bar\alpha,\bar\beta)\leq
d(\alpha,\beta)$ the projection $p\:(\hom(K,G),d)\to (\chom(K,G),\bar d)$ is continuous,
so the quotient topology on $\chom(K,G)$ is finer than the metric topology
(one can check that they are equal but we shall not use that).

The space $\chom(K,G)$ has at most countably many points \cite[Prop.\,10.14]{Ada}.
Therefore, the set $\cald=\{d(a,b)\mid a,b\in\chom(K,G)\}$ is at most
countable. Let $a,b\in\chom(K,G)$ with $a\neq b$. There exists $\lambda\in\bbr$
with $0<\lambda<\bar d(a,b)$ and $\lambda\notin\cald$.
The space $\chom(K,G)$ is then the disjoint
union of $\{x\mid \bar d(a,x)<\lambda\}$ and $\{x\mid \bar d(a,x)>\lambda\}$.
These are non-empty open sets for the topology induced by $\bar d$
and then for the quotient topology. This proves that any subspace of $\chom(K,G)$ containing more than one point is not connected.
\cqfd\sk{2}

\begin{Lemma}\label{relconju}
Let $K$ and $G$ be compact Lie groups.
Let $B$ be a space homeomorphic to a compact disk and let $b\in B$.
Let $x\mapsto\beta_x$ be a continuous map from $B$ to $\hom(K,G)$.
Then, there is a continuous  $x\mapsto g_x$ from $B$ to $G$
with $g_b=1$,
such that $\beta_x(\gamma)=g_x^\mun\,\beta_b(\gamma)\,g_x$ for all
$\gamma\in K$ and all $x\in B$.
\end{Lemma}

\preu
If $B$ is of dimension $n$, then, by a pointed homeomorphism,
one can replace the pair $(B,b)$ by $([0,1]^n,0)$ if $b$ lies in the
boundary of $B$, or by $([-1,1]^n,0)$ otherwise.
By \lemref{condiscret-bis}, $\beta_x$ stays
for all $x$ in the same
conjugacy class $\calo$ of $\hom(K,G)$. As seen in the proof of
\lemref{condiscret-bis}, the space $\hom(K,G)$ is metric, therefore
Hausdorff. Therefore, $\calo$ is compact. As $G$ is compact,
the map $p\:G\to\calo$ given by $g\mapsto g\,\beta_0\,g^\mun$
can then be identified with a principal bundle
whose structure group is the centraliser $\calz(\beta_0(K))$,
which is a closed subgroup of $G$.
\lemref{relconju} then follows from the a recursive use of
the homotopy lifting property. \cqfd

\begin{Lemma}\label{Kabelian}
Let $K$ be a compact abelian Lie group. Denote by $K_1$ the connected
component of the unit element $1\in K$. Then, there is a
bicontinuous isomorphism
$K_1\times\pi_0(K) \xrightarrow{\approx} K$.
\end{Lemma}

\preu As $K$ is abelian, it suffices to construct a homomorphic section
of the projection $K\to\pi_0(K)$. As $\pi_0(K)$ is finite, one can reduce
to the case where $\pi_0(K)=C$ is a cyclic group of order $m$.
Let $c\in C$ be a generator and choose $\tilde c\in K$ representing $c$.
Then, $\tilde c^m$ is in $K_1$ and there exists
$\gamma\in K_1$ such that $\gamma^m=\tilde c^m$. The map
$\sigma\:C\to G$ defined by $\sigma(c^k)=\tilde c^k\gamma^{-k}$ is
a homomorphic section of the projection $K\to C$. \cqfd

\begin{ccote}\emph{Spaces over $A$}. Let $f\:X\to A$ be a continuous map between
topological spaces. This enables us to consider $X$ as a ``space
over $A$''. For $a\in A$, the {\it stalk} over $a$ is $X_a=f^\mun(a)$.
Any subspace $Y$ of $B\times A$ is seen as a space over $A$ via the projection
onto $A$ restricted to $Y$.
\end{ccote}

\section{Split $\Gamma$-spaces}\label{Ssplsp}

Let $A$ be a topological space and $\Gamma$ be a topological group.
A {\it $\Gamma$-space} is a topological space equipped with a continuous
left action of $\Gamma$.
If $X$ is a $\Gamma$-space and $x\in X$, we denote by $\Gamma_x$
the stabiliser of $x$.

A {\it \spl $\Gamma$-space over $A$} is a triple $(X,\pi,\varphi)$
where
\begin{enumerate}
\item $X$ is a $\Gamma$-space.
\item $\pi\:X\to A$ is a continuous surjective map and, for each $a\in A$, 
the preimage $\pi^\mun(a)$ is a single orbit. 
\item $\varphi \: A\to X$ is a continuous section of $\pi$
\end{enumerate}

The maps $\pi$ and $\varphi$ may ommitted from the notation and we might speak
of a \spl $\Gamma$-space $X$ over $A$.
By (ii), the map $\pi$ descends to  
$\bar\pi\:\Gamma\backslash X\to A$
which is a homeomorphism (its continuous inverse
is provided by the section $\varphi$). The map
$\pi$ may thus be identified with the projection of $X$ onto
the orbit space $\Gamma\backslash X$.

A {\it $(\Gamma,A)$-groupoid} is a subspace $\cali\subset \Gamma\times A$
such that, for each $a\in A$, the space $\cali_a=\cali\cap(\Gamma\times\{a\})$
is of the form $\tilde\cali_a\times\{a\}$, 
where $\tilde\cali_a$ is a closed subgroup of $\Gamma$. We consider
$\cali_a$ as a topological group, naturally isomorphic to the closed
subgroup $\tilde\cali_a$ of $\Gamma$. We will often identify these two groups, and write, for instance, $\cali_a=\cali_b$ when we mean $\tilde\cali_a=\tilde\cali_b$.
The space $\cali$ is regarded as a topological groupoid whose space of objects is $A$: 
if $a,b\in A$, the space of morphisms from $a$ to $b$ is empty
when $a\neq b$ and is equal to $\cali_a$ otherwise.

Let $\cali$ be a $(\Gamma,A)$-groupoid.
A {\it (left) action} of $\cali$ on a topological space $W$ is a
continuous map $\beta\:\cali\times W\to W$ such that, for each $a\in A$,
the restriction of $\beta$ to $\cali_a\times W$ is an action of
$\cali_a$ on $W$. The notation $\zeta\cdot w$ is used
for $\beta(\zeta,w)$. A right action is defined accordingly,
as a continuous map from $W\times \cali$ to~$W$.
When $\cali$ acts on the right on a space $V$ and on the left
on a space $W$, we define the quotient space
$$
V\times_{\cali} W = (V\times W)/\sim \, ,
$$
where ``$\sim$'' is the smallest equivalence relation such that
$(v\cdot\zeta,w)\sim (v,\zeta\cdot w)$ for all $\zeta\in \cali$.
If $W$ is a space over $A$,
then $V\times_{\cali} W$ is a space over $A$ as well.
The stalk over $a$ is then $V\times_{\cali_a} W_a$.

Let $(X,\pi,\varphi)$ be a \spl $\Gamma$-space over $A$.
Its {\it isotropy groupoid} is the $(\Gamma,A)$-groupoid
defined by
$$
\cali(X)=\cali(X,\pi,\varphi)\:=\{(\gamma,a)\in \Gamma\times A
\mid \gamma\in\Gamma_{\varphi(a)}\} \ .
$$

A $(\Gamma,A)$-groupoid $\cali$ is called {\it \loma\ } if each point $a\in A$
admits a neighbourhood $U$ such that $\cali_u$ is a subgroup of $\cali_a$
for all $u\in U$. A space $X$ is called {\it \lcom}
if it is Hausdorff and if every point of $X$ admits a compact neighbourhood.
The main result of this section is the following 
proposition.

\begin{Proposition}[Reconstruction]\label{recon-pro}
Let $\Gamma$ be a compact topological group and $A$ be a \lcom\ space.
Let $\cali$ be a \loma\ $(\Gamma,A)$-groupoid. Then, the 
following properties hold.

\begin{enumerate}
\item There is a \spl $\Gamma$-space
$(Y_\cali,\varpi,\phi)$ over $A$ with isotropy groupoid $\cali$; the space
$Y_\cali$ is \lcom.
\item Let $(X,\pi,\varphi)$ and $(X',\pi',\varphi')$ be two
\spl $\Gamma$-spaces over $A$  with isotropy groupoid $\cali$.
Suppose that $X$ and $X'$ are \lcom. Then
there is a unique $\Gamma$-equivariant homeomorphism $F\:X\to X'$
such that $\varphi'=F\pcirc\varphi$.
\end{enumerate}
\end{Proposition}

\proref{recon-pro} permits us to speak about {\it the} \spl
$\Gamma$-space over $A$  with isotropy groupoid $\cali$
(as we speak about {\it the} real number field instead of 
{\it a} real number field). The triple
$(Y_\cali,\varpi,\phi)$ constructed for the proof of (i) 
is an explicit model for this space, 
but other models also occur naturally, as will be seen in examples.

\sk{1}
\ppreu{of \proref{recon-pro}}
The groupoid $\cali$ acts by multiplication on the right on $\Gamma$.
We let it act trivially on the left on $A$ and form the space
$$Y_\cali = \Gamma\times_\cali A .$$
The projection $\Gamma\times A\to A$ descends to a continuous
surjective map $\varpi\:Y_\cali\to A$. 
The section $\phi\:A\to Y_\cali$ is defined
by $\phi(a)=[1,a]$, where $1$ is the unit element in $\Gamma$.
The $\Gamma$-action $\beta\cdot (\gamma,x)=(\beta\gamma,x)$ on
$\Gamma\times A$ descends to a $\Gamma$-action on $Y_\cali$. 
The stalk $\varpi^\mun(a)$ is the orbit through $\phi(a)$
and $\Gamma_{\phi(a)}=\cali_a$. Thus,
$(Y_\cali,\varpi,\phi)$ is a \spl $\Gamma$-space over $A$ 
with isotropy groupoid $\cali$.

To prove that $Y_\cali$ is Hausdorff, let $x$ and $y$
be two distinct points in $Y_\cali$. Let us show that they 
admit disjoint neighbourhoods. If $\pi(x)\neq\pi(y)$,
this is obvious since $A$ is Hausdorff. In the case $\pi(x)=\pi(y)=a$,
let $U$ be a neighbourhood of $a$ for which $\cali_b$
is a subgroup of $\cali_a$ for all $b\in U$. Then, $\varpi^\mun(U)$
is a neighbourhood of $\{x,y\}$ and there is a continuous map
$f\:\varpi^\mun(U)\to (\Gamma_{\phi(a)}\backslash\Gamma)\times U$
such that $f(x)\neq f(y)$. As $\Gamma_{\phi(a)}=\cali_a$ is a
closed subgroup of $\Gamma$, the space  
$(\Gamma_{\phi(a)}\backslash\Gamma)\times U$ is Hausdorff
and then $x$ and $y$ admit disjoint neighbourhoods. Observe that
the proof that $Y_\cali$ is Hausdorff uses only that
$\Gamma$ and $A$ are Hausdorff and that $\cali$ is \loma.

Let $x\in Y_\cali$. As $A$ is \lcom, $\varpi(x)$ admits a compact 
neighbourhood $V$ in $A$. Then, $\varpi^\mun(V)$ is
a neighbourhood of $x$ and is the continuous image of $\Gamma\times V$
under the natural projection $\Gamma\times A\to Y_\cali$. 
As $\Gamma\times V$ is compact and $Y_\cali$ is Hausdorff, 
$\varpi^\mun(V)$ is compact. This shows that $Y_\cali$ is \lcom\
and finishes the proof of (i).

Let us prove (ii), starting with uniqueness. Let $F_1,F_2\:X\to X'$ be two $\Gamma$-equivariant isomorphisms
satisfying $F_1\pcirc\varphi=F_2\pcirc\varphi$. Then 
$F_1=F_2$ by equivariance, since $\varphi(A)$
is a fundamental domain for the $\Gamma$-action.
Now, let $(X,\pi,\varphi)$ be a \spl $\Gamma$-space with 
isotropy groupoid $\cali$. 
The map $\tilde F\:\Gamma\times A\to X$ defined by 
$\tilde F(\gamma,a)=\gamma\cdot\varphi(a)$ is continuous, 
$\Gamma$-equivariant and surjective.
It descends to a $\Gamma$-equivariant continuous bijection $F\:Y_\cali\to X$.
Observe that $\tilde F$ is proper, that is $\tilde F^\mun(K)$
is compact for all compact subsets $K$ in $X$. Indeed, 
$\tilde F^\mun(K)$ is a closed subset of $\Gamma\times\pi(K)$
which is compact. Since $Y_\cali$ is Hausdorff, the map 
$F$ is proper as well. A proper continuous bijection 
between \lcom\ spaces is a homeomorphism,
which proves (b). \cqfd 

\begin{Remark}
If, in \proref{recon-pro}, we only assume that $\Gamma$ is locally
compact, the space $Y_\cali$ constructed for the proof of (i) might
not be \lcom. As an example, take 
$A=[0,1]$, $\cali_a$ trivial for $a<1$ and $\cali_1=\Gamma$.
Then $Y_\cali$ is the cone on $\Gamma$, which is not \lcom\ if 
$\Gamma$ is not compact. 
Moreover, the uniqueness  also fails in this case.
Let $\Gamma=\bfz$. Then the cone $V$ on the real integers  in
the complex plane, with vertex $i$ say, and the induced metric, is a split
$\Gamma$-space with isotropy groupoid $\cali$.
The proof of Proposition 2.1 provides a $\Gamma$-equivariant continuous
bijection from the $Y_\cali$ onto $V$ but it is not a homeomorphism
(a set containing one point in the interior of each segment
would be closed in $Y_\cali$, even if it contains a subsequence converging to
the vertex $i$). A version of \proref{recon-pro} with $\Gamma$ non-compact
is given in \proref{recon-pro-CW}.
\end{Remark}

A stronger local maximality condition will play a role
in Section~\ref{nspeqbd}.
A $(\Gamma,A)$-groupoid $\cali$ is called {\it \sloma}\ if,
for each point $a\in A$ and each neighbourhood $B$ of $a$,
there exists an open set $U$ of $A$, with $a\in U\subset B$, 
and a homotopy $\rho_t\:U\to U$ ($t\in [0,1]$) satisfying
$\rho_0(u)=u$, $\rho_1(u)=a$ and $\cali_u\subset \cali_{\rho_t(u)}$
for all $u\in U$ and $t\in [0,1]$. The neighbourhood $U$ is then contractible.
\Sloma\ implies \loma.

\begin{Lemma}\label{splonsloma}
Let $\Gamma$ be a compact topological group and $A$ be a \lcom\ space.
Let $\cali$ be a \sloma\ $(\Gamma,A)$-groupoid.
Let $(X,\pi,\varphi)$ be a \spl $\Gamma$-space over $A$
with isotropy groupoid $\cali$ and let $x\in X$. Then, there is
a $\Gamma$-equivariant open neighbourhood $\hat U$ of the orbit $\Gamma x$
and a $\Gamma$-equivariant $\hat\rho_t\:\hat U\to\hat U$ such that
$\hat\rho_0={\rm id}$ and $\hat\rho_1(\hat U)=\Gamma\cdot x$.
Moreover, $\hat\rho_t$ satisfies $\pi\pcirc\hat\rho_t=\rho_t\pcirc\pi$
and  $\varphi\pcirc\rho_t=\hat\rho_t\pcirc\varphi$.
\end{Lemma}

\preu By \proref{recon-pro}, one may suppose that
$(X,\pi,\varphi)=(Y_\cali,\varpi,\phi)$.
Let $a=\varpi(x)$ and let $\rho_t\:U\to U$ be a homotopy from an open
neighbourhood $U$ of $A$ to itself, such that
$\rho_0(u)=u$, $\rho_1(u)=a$ and $\cali_u\subset \cali_{\rho_t(u)}$
for all $u\in U$. Let $\hat U=\varpi^\mun(U)$. We check that
the required homotopy $\hat\rho_t$ can be defined by
$\hat\rho_t([\gamma,u])\:=[\gamma,\rho_t(u)]$. \cqfd 

Let $\cali$ be a $(\Gamma,A)$-groupoid and let $\cali'$ be a 
$(\Gamma',A')$-groupoid. A {\it morphism of groupoids}
from $\cali$ to $\cali'$ is a commutative diagram 
$$
\xymatrix@C-3pt@M+2pt@R-4pt{%
\cali \ar[d]
\ar[r]^(0.50){f}  &
\cali' \ar[d]  \\
A \ar[r]^{\bar f}  &
A'
}
$$
where $f$ and $\bar f$ are continuous maps, such that,
for each $a\in A$, the restriction of $f$ to $\cali_a$
is a homomorphism $f_a\:\cali_a\to\cali'_{\bar f(a)}$.
The map $\bar f$ is not mentioned when it is obvious, like
an inclusion or a constant map.

\section{Split equivariant principal bundles}\label{nspeqbd}

\subsection{Definitions}\label{nspeqbd-defi}
Let $G$ be a topological group and $X$ be a topological space.
By a {\it $G$-principal bundle $\eta$ over $X$},
we mean, as usual, a continuous surjection $p\colon  E\to X$ from
a space $E=E(\eta)$ and a free right action $E\times G\to E$
so that $p(z\cdot g) = p(z)$, with the standard
local triviality condition.
Two $G$-principal bundles $\eta\: E\xrightarrow{p} X$ and
$\eta'\: E'\xrightarrow{p'} X$ over $X$
are {\it isomorphic} if there exists a $G$-homeomorphism
$f\:E\to E'$ such that $p'\pcirc f= p$. Isomorphism classes
of $G$-principal bundles over $X$ are denoted by $\bung(X)$.

Let $X$ be a $\Gamma$-space for a topological group $\Gamma$.
A $G$-principal bundle  $\eta\: E\xrightarrow{p} X$ is 
called a {\it $\Gamma$-equivariant principal $G$-bundle}
if it is given a left action $\Gamma\times E\to E$
commuting with the free right action of $G$ and
such that the projection $p$ is $\Gamma$-equivariant.
Two $\Gamma$-equivariant principal $G$-bundles $\eta$
and $\eta'$ are called {\it $\Gamma$-isomorphic} 
(or just {\it isomorphic}) if there exists
a  $G$-homeomorphism from $E(\eta)$ to $E(\eta')$ 
over the identity of $X$ which is $\Gamma$-equivariant.
The set of $\Gamma$-isomorphism classes of $\Gamma$-equivariant
$G$-principal bundles over $X$ is denoted by $\nbungt(X)$.
There is a forgetful map $\nbungt(X)\to\bung(X)$.

Let $(X,\pi,\varphi)$ is a \spl $\Gamma$-space over $A$
with isotropy groupoid $\cali$.
Let $\xi$ be a $\Gamma$-equivariant principal $G$-bundle over $X$.
We say that $\xi$ is {\it \spl} if
the induced bundle $\varphi^*\xi$ is trivial.
For instance, any $\Gamma$-equivariant principal $G$-bundle is \spl when
$A$ is contractible and paracompact, which is the case in many examples
of Subsection~\ref{Spltsp-ex}.
Two \spl $\Gamma$-equivariant principal $G$-bundles 
over $(X,\pi,\varphi)$ are {\it isomorphic} if they 
are isomorphic just as $\Gamma$-equivariant principal $G$-bundles
over $X$.
The set of isomorphism classes of $\Gamma$-equivariant \spl
$G$-principal bundles over $(X,\pi,\varphi)$ is denoted by
$\nsbungt(X,\pi,\varphi)$ or simply by $\nsbungt(X)$.
It is a subset of $\nbungt(X)$.

\subsection{The isotropy representation}\label{S-isor}

Let $\cali$ be a $(\Gamma,A)$-groupoid and $G$ be a topological group.
A {\it continuous representation} of $\cali$ to $G$
is a continuous map $\alpha\:\cali\to G$ such that, for all $a\in A$,
the restriction $\alpha_a$ of $\alpha$ to $\cali_a$ is a homomorphism
(it is thus a morphism of groupoids between $\cali$ and the 
$(pt,G)$-groupoid $G\to pt$).
Two continuous representations
$\alpha_1$ and $\alpha_2$ are called {\it conjugate} if there
exists a continuous map $\psi\:A\to G$
such that $\alpha_2(\zeta)=\psi(\pi_2(\zeta))^\mun\alpha_1(\zeta)\psi(\pi_2(\zeta))$
for all $\zeta\in \cali$, where $\pi_2\colon \cali \to A$ is the second factor projection.

A continuous representation of $\alpha\colon \cali \to G$
is called {\it \rloma}\ if each point $a\in A$
admits a neighbourhood $U$ such that $\cali_u$ is a subgroup of $\cali_a$
for all $u\in U$, together with a continuous map $g\:U\to G$ such that
$\alpha_u(\gamma)=g(u)\alpha_a(\gamma)g(u)^\mun$ for all
$u\in U$ and all $\gamma\in\cali_u$. This implies that
$\cali$ is \loma. It is easy to see that,
if $\alpha,\beta\:\cali\to G$ are two conjugate representations
of $\cali$, then $\beta$ is \rloma\ if and only if $\alpha$ is \rloma.
We denote by $\repi$ the set of conjugacy classes of
\rloma\ continuous representations of $\cali$.

Let $(X,\pi,\varphi)$ be a \spl $\Gamma$-space over $A$
with isotropy groupoid $\cali$.
Let $\eta\: E\xrightarrow{p} X$ be a \spl
$\Gamma$-equivariant $G$-principal bundle over $X$.
As $\varphi^*\eta$ is trivial, there exists a continuous lifting
$\tilde\varphi\:A\to E$ of $\varphi$.
The equation
\begin{equation}\label{conjueq}
\gamma\cdot\tilde\varphi(a) = \tilde\varphi(a)\,\tilde\alpha_a(\gamma) \, ,
\end{equation}
valid for $a\in A$ and $\gamma\in\cali_a$, determines a continuous
representation $\alpha_{\eta,\tilde\varphi}\:\cali\to G$.

\begin{Lemma}\label{rloma}
Suppose that $\Gamma$ and $G$ are compact Lie groups and that $A$ is \lcom.
If $\cali$ is \sloma, then the continuous representation
$\alpha_{\eta,\tilde\varphi}$ is \rloma.
\end{Lemma}

\preu Let $a\in A$ and let $B$ be a compact neighbourhood of $a$.
Since $\cali$ is \sloma, there exists an open set $U$,
with $a\in U\subset B$ and a homotopy
$\rho_t\:U\to U$  such that
$\rho_0(u)=u$, $\rho_1(u)=a$ and $\cali_u\subset \cali_{\rho_t(u)}$
for all $u\in U$. If $Z\subset A$, we denote $\hat Z=\pi^\mun(Z)$;
if $Y$ is a $\Gamma$-invariant subspace of $X$, we denote
$E_{Y}=p^\mun(Y)$. The latter
is the total space of a \spl$\Gamma$-equivariant $G$-principal
bundle $\eta_{\scr Y}$ over $Y$.

By \proref{recon-pro} and its proof, the space $\hat B$ is compact.
Then, $E_{\hat B}$ is compact and therefore totally regular.
By \cite[Proposition~8.10]{D2},
the bundle $\eta_{\scr\hat B}$ is then a locally trivial numerable
$\Gamma$-equivariant $G$-principal bundle
in the sense of \cite[p.~58]{D2}. The same then holds for
its restriction $\eta_{\scr\hat U}$.

By \lemref{splonsloma} and its proof, the homotopy
$\rho_t\:U\to U$ is covered by a $\Gamma$-equivariant
homotopy $\hat\rho_r\:\hat U\to\hat U$ such that
$\hat\rho_0={\rm id}$ and $\hat\rho_1(\hat U)=\pi^\mun(a)$.
By \cite[Theorem~8.15]{D2}, the induced bundle
$\hat\rho_1^*\eta_{\pi^\mun(a)}$ is then isomorphic
to $\eta_{\hat U}$. More precisely, let
$$
E_1\:=\{(x,z)\in\hat U\times E_{\pi^\mun(a)}\mid \hat\rho_1(x)=p(z)\}
$$ 
be the total space of $\hat\rho_1^*\eta_{\pi^\mun(a)}$. Then,
there is a $(\Gamma\times G)$-equivariant homeomorphism
$\mu\: E_1\to E_{\hat U}$ which commutes with the projections onto $\hat U$.
By \lemref{splonsloma} one has $\varphi\pcirc\rho_t=\hat\rho_t\pcirc\varphi$;
therefore, $(\varphi(u),\varphi(a))\in E_1$ for all $u\in U$.
This enables to define $\tilde\varphi'\:U\to E$ by
$\tilde\varphi'(u)=\mu(\varphi(u),\tilde\varphi(a))$.
For $\gamma\in\cali_u\subset\cali_a$, we have
\begin{equation}\label{rloma-eq1}
\gamma\cdot\tilde\varphi'(u)=\mu(\varphi(u),\gamma\tilde\varphi(a))
=\mu(\varphi(u),\tilde\varphi(a)\alpha_a(\gamma))=
\mu(\varphi(u),\tilde\varphi(a))\cdot\alpha_a(\gamma) \, .
\end{equation}
On the other hand, $\tilde\varphi$ and $\tilde\varphi'$ are two
liftings of $\varphi$ over $U$. Hence, there exists a continuous
map $g\:U\to G$ such that $\tilde\varphi'(u)=\tilde\varphi(u)\cdot g(u)$
for all $u\in U$. Therefore
\begin{equation}\label{rloma-eq2}
\gamma\cdot\tilde\varphi'(u)=\gamma\tilde\varphi(u)\,g(u)=
\tilde\varphi(u)\alpha_u(\gamma)g(u)=
\tilde\varphi'(u)\cdot\big(g(u)^\mun\alpha_u(\gamma)g(u)\big) \, .
\end{equation}
Comparing Equations~\eqref{rloma-eq1} with~\eqref{rloma-eq2},
we get that $\alpha_u(\gamma)=g(u)\alpha_a(\gamma)g(u)^\mun$
which proves \lemref{rloma}. \cqfd

By \lemref{rloma}, $\alpha_{\eta,\tilde\varphi}$ determines a class
$\alpha_\eta\in\repi$. We check that $\alpha_{\eta}$
does not depend on the choice of $\tilde\varphi$
and depends only on the $\Gamma$-equivariant isomorphism class of $\eta$; 
details are as in the proof of \cite[Lemma~3.2]{hhausmann2} .
This defines a map $$\isor\:\nsbungt(X)\to\repi$$ called the 
{\it isotropy representation}.

\subsection{The classification theorem}\label{S-classi}

The following theorem corresponds to Theorem~A of the introduction.

\begin{Theorem}[Classification]\label{Th-classi}
Let $(X,\pi,\varphi)$ be a \spl $\Gamma$-space over $A$ with isotropy groupoid $\cali$.
Suppose that $A$ is \lcom, that $\cali$ is \sloma\ and that
$\Gamma$ is a compact Lie group. Then, for any compact Lie group $G$, 
the isotropy representation $\isor\:\nsbungt(X)\to\repi$ is a bijection.
\end{Theorem}

\preu
We first prove the surjectivity of $\isor$. By \proref{recon-pro}, we may assume that
$(X,\pi,\varphi)=(Y_\cali,\varpi,\phi)$. Recall that $Y_\cali=\Gamma\times_{\cali}A$.

Let $\beta\:\cali\to G$ be a continuous representation.
Then $\cali$ acts on the left on $A\times G$ by $\zeta\cdot (a,g)=(a,\beta(\zeta)\, g)$.
Form the space $E_\beta=\Gamma\times_{\cali}(A\times G)$. The continuous map
$p\:E_\beta\to Y_\cali$ given $p([\gamma,(a,g)])=[\gamma,a]$ coincides
with the projection $E_\beta\to E_\beta/G$ of $E_\beta$ to its orbit
space for the obvious free right $G$-action on $E_\beta$.
A lifting $\tilde\phi\:A\to E_\beta$ of $\phi$ is given by
$\tilde\phi(a)=[1,(a,1)]$, where $1$ denotes the unit elements.
For $a\in A$ and $\gamma\in \cali_a$, one has
$$
\gamma\cdot\tilde\phi(a) =\gamma\cdot[1,(a,1)]=[\gamma,(a,1)]=[1,(a,\beta_a(\gamma)]
=\tilde\phi(a)\cdot \beta_a(\gamma)\ .
$$
We now prove that $p$ admits local trivializations when $\beta$ is \sloma.
Let $a\in A$. Choose an open neighbourhood $U_a$ of $a$ 
such that $\cali_u$ is a subgroup of $\cali_a$
for all $u\in U_a$, together with a continuous map $g_a\:U_a\to G$ such that
$\beta_u(\gamma)=g_a(u)\beta_a(\gamma)g_a(u)^\mun$ for all 
$u\in U_a$ and all $\gamma\in\cali_u$. 
This gives an open cover $\calu=\{U_a\mid a\in A\}$ of $A$.
Setting $\hat U_a=\varpi^\mun(U_a)$
gives rise to an open cover $\hat\calu=\{\hat U_a\mid a\in A\}$ of 
$X$, indexed by $A$.
Define $\tilde f_a\:\Gamma\times U_a \times G \to \Gamma\times \{a\} \times G$ by
$f_a(\gamma,u,g)\:=(\gamma,a,g_a(u)g)$. If $\delta\in\cali_u$, we have
$$
\tilde f_a(\gamma\delta,u,g)=(\gamma\delta,a,g_a(u)g)=(\gamma,a,\beta_a(\delta)g_a(u)g)
$$
and
$$
\tilde f_a(\gamma,u,\beta_u(\delta)g)=(\gamma,a,g_a(u)\beta_u(\delta)g)\ .
$$
Since $\beta_a(\delta)g_a(u)= g_a(u)\beta_u(\delta)$, this shows that $\tilde f_a$
descends to a continuous $G$-equivariant map
$f_a\:p^\mun(\hat U_a) \to p^\mun(\pi^\mun(a))$. Passing to the quotient by $G$ gives
rise to a commutative diagram
$$\xymatrix{p^\mun(\hat U_a)\ar[r]^{f_a}\ar[d]^p&p^\mun(\pi^\mun(a))
\ar[d]^p& \Gamma\times_{\cali_a}G\ar[l]^(0.4){\approx}\ar[d]\\
\hat U_a\ar[r]^{\bar f_a}&\ \pi^\mun(a)\ &\Gamma/\cali_a\ar[l]^(0.4){\approx}
}
$$
Since $\Gamma$ is a Lie group and $\cali_a$ a closed subgroup, 
the projection $q_a\:\Gamma\to\Gamma/\cali_a$ admits local sections
$\sigma_{\scr V}\:V\to\Gamma$ for each $V$ in some open covering
$\calv_a$ of $\Gamma/\cali_a$ (see, e.g. \cite[\S\,7.5]{St}).
We check that the formula
$$
\zeta_{\scr V}(\gamma,g)=\beta_a(\sigma(q(\gamma))^\mun\gamma)\,g
$$
defines a $G$-equivariant continuous map $\zeta_{\scr V}\:p^\mun(V)\to G$,
which gives rise to a trivialization over $V$ of
$p\:\Gamma\times_{\cali_a}G\to \Gamma/\cali_a$.
Therefore, $\zeta_{\scr V}\pcirc f_a\:p^\mun(\bar f_a^\mun(V))\to G$
is a $G$-equivariant continuous map giving rise to
a trivialization over $\bar f_a^\mun(V)$ of $p\: p^\mun(\hat U_a)\to \hat U_a$.
This gives rise to a trivializing open cover 
$\calw=\{\bar f_a^\mun(V)\mid (a,V)\in \cala \}$ of $X$, indexed by
$\cala=\{(b,V)\mid b\in A \hbox{ and } V\in\calv_b\}$.
We have proved thus the surjectivity of $\isor$.

We now prove the injectivity of $\isor$.
Let $\eta\:(E\xrightarrow{\bar p} X)$ be a \spl $\Gamma$-equivariant
bundle with $\isor(\eta)=[\beta]$. A lifting $\bar\varphi\:A\to E$ of
$\varphi$ then produces a continuous representation $\bar\beta=\alpha_{\eta,\bar\varphi}$
with $[\bar\beta]=[\beta]$. There exists then a continuous map $\psi\:A\to G$
such that $\beta(\zeta)=\psi(q(\zeta))^\mun\bar\beta(\zeta)\psi(q(\zeta))$.
The map $\tilde\varphi\:A\to E$ given by $\tilde\varphi(a)=\bar\varphi(a)\cdot\psi(a)$
is then another lifting of $\varphi$ such that
$\alpha_{\eta,\tilde\varphi}=\beta$.
One checks that the correspondence $[\gamma,(a,g)]\mapsto \gamma\cdot\varphi(a)\cdot g$
defines a $(\Gamma\times G)$-equivariant continuous bijection $\tilde F\:E_\beta\to E$, covering the unique $\Gamma$-equivariant homeomorphism $F\:Y_\cali\to X$ such that
$F\pcirc\varphi=\phi$, obtained in \proref{recon-pro}.
Since $F$ is a homeomorphism, so is $\tilde F$. Indeed, choose an open set
$Z$ in $Y_\cali$ such that $\xi_\beta$ is trivial over $Z$ and
$\xi$ is trivial over $F(Z)$. Using trivializations, we can write
$\tilde F(z,g)=(F(z),\mu(z)\,g)$, where $\mu\:Z\to G$ is a continuous map.
Then $\tilde F^\mun$ has, over $F(Z)$, the form
$\tilde F^\mun(y,h)=(F^\mun(y),\mu(F^\mun(y))^\mun\,h)$ which is continuous.
We have thus proven that two \spl $\Gamma$-equivariant
principal $G$-bundles $\eta$ and $\eta'$ with $\isor(\eta)=\isor(\eta')$
are $\Gamma$-equivariantly isomorphic. \cqfd

\begin{Remark}\label{numerable} 
Recall that an open cover of a space is {\it numerable} if it admits a refinement
by a locally finite partition of unity. In the proof of \thref{Th-classi},
the covers $\calv_a$ are numerable, since $\Gamma/\cali_a$ are manifolds.
Hence we can check that the trivializing cover 
$\calw$ of $X$ is numerable if $\calu$ is numerable. This observation will be used in Theorem \ref{Th-classi-CW}.
\end{Remark}

\begin{ccote}\label{as-ab} \emph{Non-split bundles and abelian structure group}.
Let $(X,\pi,\varphi)$ be a \spl $\Gamma$-space over $A$. 
One has the map $\varphi^*\:\nbungt(X)\to\bung(A)$, sending $\xi$ to $\varphi^*\xi$.
This map is surjective: if $\eta\in\bung(A)$, then 
$\pi^*\eta$ admits a natural $\Gamma$-action, since $\pi$ is $\Gamma$-invariant,
so $\pi^*$ is a section of $\varphi^*$. \thref{Th-classi} computes the pre-image
of the trivial bundle, which is $\nsbungt(X)$.

Let us now assume that $G$ is abelian. Recall that there is 
then a composition law
``$\otimes$'' on  $\nbungt(X)$ which makes the latter an abelian group. 
If $\xi_i\:E_i\xrightarrow{p_i} X$ ($i=1,2$) are $\Gamma$-equivariant
principal $G$-bundles, one defines $\xi_1\otimes\xi_2\:E\xrightarrow{p_i} X$
by first forming the pull-back
$$
\xymatrix@C-3pt@M+2pt@R-4pt{%
E_1\hat\times E_2 \ar[d]
\ar[r]  &
E_1 \ar[d]^(0.50){p_1}  \\
E_2 \ar[r]^(0.50){p_2} &
X}
$$ 
where the map $E_1\hat\times E_2 \to X$ is a principal $G\times G$-bundle. Set $E=E_1\hat\times_G E_2$ (as $G$ is abelian, it acts on the left or on the right
on $E_i$) and check that $\xi_1\otimes\xi_2$ is a principal $G$-bundle over $X$.
The diagonal $\Gamma$-action on $E_1\hat\times E_2$ descends to a $\Gamma$-action 
on $E$, making  $\xi_1\otimes\xi_2$ a $\Gamma$-equivariant principal 
$G$-bundle. When $G=S^1$, we can think of $\xi_i$ as $\Gamma$-equivariant
complex line bundles over $X$, thus ``$\otimes$'' 
becomes the standard tensor product. The map $\varphi^*\:\nbungt(X)\to\bung(A)$
is a group-homomorphism.

Another special feature of the case $G$ abelian is that the isotropy representation
is defined on $\nbungt(X)$: in Equation~\eqref{conjueq}, one can just use
a local section $\tilde\varphi$ around $a\in A$, whose choice is irrelevant 
if $G$ is abelian. The set $\repi$ is an abelian group, by multiplication
of the images, and $\isor\:\nbungt(X)\to\repi$ is a group homomorphism.
Using \thref{Th-classi}, we get
\begin{Proposition}\label{Th-classi-ab}
Let $(X,\pi,\varphi)$ be a \spl $\Gamma$-space over $A$ with isotropy groupoid $\cali$.
Suppose that $A$ is \lcom, that $\cali$ is \sloma\ and that
$\Gamma$ is a compact Lie group. Then, for any compact abelian Lie group $G$, 
one has an isomorphism of abelian groups
$$
(\isor,\varphi^*)\: \nbungt(X) \xrightarrow{\approx}
\repi\times\bung(A) \, .  
\mancqfd
$$
\end{Proposition}
\end{ccote}

\begin{ccote}\label{funct} \emph{Functorial properties}.
\thref{Th-classi} enjoys functorial properties which are contravariant
in $(\Gamma,A)$ and covariant in $G$. For the contravariant ones,
let $f\:A'\to A$ be a continuous map between \lcom\ spaces and $h\:\Gamma'\to\Gamma$
be a continuous homomorphism between compact Lie groups. Let $\cali$ be a \gagrou. Then
$$
\cali':=(h,f)^*\cali := \{(\gamma',a')\in\Gamma'\times A' \mid h(\gamma')\in\cali_{f(a')}\}
$$
is a $(\Gamma',A')$-groupoid, with $\cali'_{a'}=h^\mun(\tilde\cali_{f(a')})\times\{a'\}$.
One has the continuous map
\begin{equation}\label{funct-eq1}
\Gamma'\times_{\cali'}A' \xrightarrow{(h,f)} \Gamma\times_{\cali}A
\end{equation}
Therefore, if $\cali$ and $\cali'$ are \sloma, \proref{recon-pro} 
together with Equation~\eqref{funct-eq1}
implies the following: if $(X,\pi,\varphi)$ and $(X',\pi',\varphi')$ are the split spaces with
isotropy groupoids $\cali$ and $\cali'$, there is a unique map $F=F_{h,f}\:X'\to X$ such that
$F(\gamma'x')=h(\gamma')F(x')$, $\pi\pcirc F=f\pcirc\pi'$ and $\varphi\pcirc f=F\pcirc\varphi'$.
Let $\eta$ be a \spl $\Gamma$-equivariant principal $G$-bundle over $X$.
Using \thref{Th-classi}, one checks that $\eta':=F^*\eta$ is a
\spl $\Gamma'$-equivariant principal $G$-bundle over $X'$ and that the 
isotropy representations $\alpha'\in{\rm Rep}^G(\cali')$ 
and $\alpha\in\repi$ satisfy $\alpha'=h^*\alpha$,
where $h^*\alpha=\alpha\pcirc h$.
Therefore, one gets a commutative diagram
\begin{equation}\label{funct-contradiag}
\begin{array}{c}{\xymatrix@C-3pt@M+2pt@R-4pt{%
\nsbungt(X) \ar[d]^(0.50){\approx}_(0.50){\isor}
\ar[r]^(0.50){F^*}  &
{\rm SBun}^G_{\Gamma'}(X') \ar[d]^(0.50){\approx}_(0.50){\isor}  \\
\repi \ar[r]^(0.50){h^*}  &
{\rm Rep}^G(\cali')
}}\end{array} \ .
\end{equation}

As for the covariant functoriality in $G$, let $\mu\:G\to G'$ be a continuous homomorphism
between compact Lie groups. If $\eta\:(E\to X)$ is a \spl $\Gamma$-equivariant 
principal $G$-bundle over $X$, one checks that $\mu_*\eta\:(E\times_\mu G'\to X)$
is a \spl $\Gamma$-equivariant principal $G$-bundle with isotropy representation
$\mu_*\alpha=\mu\pcirc\alpha$. One 
gets a commutative diagram
\begin{equation}\label{funct-covadiag}
\begin{array}{c}{\xymatrix@C-3pt@M+2pt@R-4pt{%
\nsbungt(X) \ar[d]^(0.50){\approx}_(0.50){\isor}
\ar[r]^(0.50){\mu_*}  &
{\rm SBun}^{G'}_{\Gamma}(X) \ar[d]^(0.50){\approx}_(0.50){\isor}  \\
\repi \ar[r]^(0.50){\mu_*}  &
{\rm Rep}^{G'}(\cali)
}}\end{array} \ .
\end{equation}
In particular, let $G=G'\times G''$ and let $p'$ and $p''$ be the two projections.
Diagram~\eqref{funct-covadiag} becomes
\begin{equation}\label{proddiag}
\begin{array}{c}{\xymatrix@C-3pt@M+2pt@R-4pt{%
\nsbungt(X) \ar[d]^(0.50){\approx}_(0.50){\isor}
\ar[r]^(0.35){(\mu'_*,\mu''_*)}  &
{\rm SBun}^{G'}_{\Gamma}(X) \times {\rm SBun}^{G''}_{\Gamma}(X)
\ar[d]^(0.50){\approx}_(0.50){\isor\times\isor}  \\
\repi \ar[r]^(0.40){(\mu'_*,\mu''_*)}_(0.40){\approx}  &
{\rm Rep}^{G'}(\cali)\times {\rm Rep}^{G'}(\cali)
}}\end{array} \ .
\end{equation}
Diagram~\eqref{proddiag} then shows that the map 
\begin{equation}\label{prod-eq}
(\mu'_*,\mu''_*)\:\,{\rm SBun}^{G'\times G''}_{\Gamma}(X)\xrightarrow{\approx}
{\rm SBun}^{G'}_{\Gamma}(X) \times {\rm SBun}^{G''}_{\Gamma}(X)
\end{equation}
is a bijection.
\end{ccote}

\section{Cellular groupoids - Examples}

\subsection{Cellular groupoids}\label{stru}
Let $A$ be a CW-complex, filtered by its skeleta $A^{(n)}$.
We denote by $\Lambda=\Lambda(A)$ the set of cells of $A$.
The dimension of a cell $e\in\Lambda$ is denoted by $d(e)$
and we set $\Lambda_n=\{e\in\Lambda\mid d(e)=n\}$. For each
$e\in\Lambda$, there exists a characteristic map
$\sigma_e\:(\bbd^{d(e)},\bbs^{d(e)-1})\to (A^{(d(e))},A^{(d(e)-1)})$, 
and $\sigma_e$ restricted to the interior of $\bbd^{d(e)}$
is an embedding whose image is denoted by $|e|$.  
For $a\in A$, we denote by $e(a)\in\Lambda$ the cell $e$ of
smallest dimension such that $a\in\sigma(e)$. 
The set $\Lambda$
is partially ordered: $f'\leq f$ if $f'$ is a {\it face} of $f$, which means that
there exists $x\in\bbs^{d(f)-1}$ such that $e(\sigma_{f'}(x))=f$.

Let $\Gamma$ be a topological group and $A$ be a CW-complex. 
A $(\Gamma,A)$-groupoid $\cali$ is called {\it cellular} if
it is \sloma\ and if $\tilde\cali_a=\tilde\cali_b$ when $e(a)=e(b)$.
We write $\cali^{(n)}$ for the restriction of $\cali$ over $A^{(n)}$.
Recall that $\cali_a=\tilde\cali_a\times\{a\}$ where 
$\tilde\cali_a\in {\rm Gr}(\Gamma)$, the poset of closed subgroups of $\Gamma$.  
One can then define a map $\tilde\cali\:\Lambda(A)\to {\rm Gr}(\Gamma)$ by
$\tilde\cali(e)=\tilde\cali_a$ for $a$ with $e(a)=e$. 
The local maximality of $\cali$ implies
that $\tilde\cali(e)\subset\tilde\cali(f)$ when $f\leq e$. Thus,
$\tilde\cali$ is a contravariant
functor from the poset $\Lambda(A)$ to the poset ${\rm Gr}(\Gamma)$. 
A cellular groupoid is a combinatorial construction.

\begin{Lemma}\label{Lcegro-fun}
The correspondence $\cali\to\tilde\cali$ is a bijection between 
the set of cellular groupoids whose object-space is $A$
and the set of contravariant functors from $\Lambda(A)$ to ${\rm Gr}(\Gamma)$.
\end{Lemma}

\preu The correspondence is clearly injective. For the surjectivity, let 
$\calf\:e\mapsto\calf_e$ be a contravariant functor from 
$\Lambda(A)$ to ${\rm Gr}(\Gamma)$. By induction on $n$, we shall construct
$\cali^{(n)}$, giving rise to a $(\Gamma,A)$-groupoid $\cali$, with
$\tilde\cali=\calf$ and then check that $\cali$ is \sloma. 
Define $\cali^{(0)}=\coprod_{v\in\Lambda_0(A)}\calf_v\times\{v\}$.
Suppose that $\cali^{(n-1)}$ is constructed. The $n$-skeleton $A^{(n)}$ of $A$ 
is obtained as the quotient space
$$
A^{(n)}=\big(\hskip-6pt\coprod_{e\in\Lambda_n(A)} \bbd_e\, \big) \,\, {\textstyle\coprod}\,\,
A^{(n-1)}\bigg/\{x\sim \sigma_e(x)\mid x\in\bbs_e\}
$$
where $(\bbd_e,\bbs_e)$ is a copy of $(\bbd^n,\bbs^{n-1})$ and 
$\sigma_e\:\bbd_e\to A$ is a characteristic map for the cell $e$.
We then define 
$$
\cali^{(n)}=\big(\hskip-6pt\coprod_{e\in\Lambda_n(A)}(\calf_e\times \bbd_e )\, \big )\,\, {\textstyle\coprod}\,\,
\cali^{(n-1)}\bigg/\{(\gamma,x)\sim (\gamma,\sigma_e(x))\mid x\in\bbs_e\}  \ .
$$
The equivalence relation $\sim$ makes sense since, for $x\in\bbs_e$, one has
$\calf_e\subset\tilde\cali^{(n-1)}_{\mu_e(x)}$. Clearly, $\tilde\cali=\calf$.
Now, each $a\in A$
admits a fundamental system of open neighbourhoods $U$ of $a$ such that
$e(a)\leq e(u)$ for all $u\in U$. 
One can also require that $U$ admits a
homotopy  $\rho_t\:U\to U$ such that $\rho_0={\rm id}$,
$\rho_1(U)=\{a\}$ and $e(\rho_t(u))=e(u)$ for $t<1$ 
(see \cite[Theorem~6.1 and its proof]{LW}, or proof of \lemref{rloma-cell} below).
Therefore, $\cali$ is \sloma. \cqfd

The notation $\tilde\cali$ was introduced in order to state and prove \lemref{Lcegro-fun}
properly. In future occurrences, we shall write $\cali(e)$ instead of $\tilde\cali(e)$.

\begin{Remark}\label{gammaCW}
Let $(X,\pi,\varphi)$ be a \spl $\Gamma$-space over a CW-complex
$A$, with a cellular isotropy groupoid $\ciat$. Then,
$X$ is provided with a $\Gamma$-equivariant CW-complex structure (see, e.g.
\cite[Chapter~2]{D2})
with $\Gamma$-cells indexed by $\Lambda(A)$. If $\sigma_e\:\bbd^{d(e)}\to A$ is a
characteristic map for $e\in\Lambda(A)$, then $\tilde\sigma_e\:\Gamma/\cali(e)\times \bbd^{d(e)}\to X$, defined by $\tilde\sigma(\gamma,a)=\gamma\,\varphi(a)$,
is a characteristic map for the $\Gamma$-cell of $X$ corresponding to $e$.

On the other hand, let $X$ be a $\Gamma$-CW-complex and $A=\Gamma/X$ be its orbit space with
the induced CW-structure. Suppose that there exists a section $\varphi\:A\to X$ of the 
projection $\pi\:X\to A$, so that the isotropy groupoid $\cali$ is \loma. Then $\cali$
is cellular, since  $\cali_a$ is constant on the interior of each cell. We call $X$ a \emph{\spl $\Gamma$-CW-complex} over $A$.

However, one has the following example of a \spl\ $\Gamma$-space over a CW-complex
admitting no splitting for which the isotropy groupoid is cellular.
Let $X=([0,1]\times S^2)/\{(1,x)\sim (0,-x)\}$,
the mapping cylinder of the antipodal map of $S^2$, endowed with the natural action of $\Gamma=SO(3)$.
Then $A=[0,1]/\{0\sim 1\}\approx S^1$. Any splitting is of the form
$\varphi(t)=(t,f(t))$ with $\lim_{t\to 0}f(t)=-\lim_{t\to 1}f(t)$. Thus,
$f(t)$ is not constant and $\cali$ is not \loma. Observe that $X$ is a smooth closed $3$-manifold
and that the $SO(3)$-action is smooth with cohomogeneity one.
\end{Remark}

For cellular groupoids we have a stronger version of 
\proref{recon-pro} which applies to any topological group $\Gamma$.

\begin{Proposition}[Reconstruction II]\label{recon-pro-CW} 
Let $\Gamma$ be a topological group, and $A$ be a CW-complex. 
Given a cellular $(\Gamma, A)$-groupoid $\cali$, there is a unique \spl $\Gamma$-CW-complex over $A$
with isotropy groupoid $\cali$.
\end{Proposition}
\begin{proof}
The space $Y_\cali$ is a \spl $\Gamma$-CW-complex over $A$. Suppose that
$(X,\pi,\varphi)$ is another \spl $\Gamma$-CW-complex
with isotropy groupoid $\cali$. As in the proof of \proref{recon-pro},
the map $\tilde F\colon \Gamma \times A \to X$ defined by $\tilde F(\gamma, a) = \gamma\cdot\varphi(a)$ descends to give a continuous $\Gamma$-equivariant
bijection $F\colon Y_\cali \to X$. For each cell $e\in\Lambda_n(Y)$ with characteristic map $\sigma_e\:\bbd^{n}\to A$, there is a commutative diagram
$$
\begin{array}{c}{\xymatrix@C-3pt@M+2pt@R-4pt{%
\Gamma\times_{\cali(e)} D^n \ar[d]^(0.50){\sigma_e^{Y_\cali}}
\ar[r]^(0.50){{\rm id}}  &
\Gamma\times_{\cali(e)} D^n \ar[d]^(0.50){\sigma_e^Y}  \\
X \ar[r]^(0.50){F}  &
Y \, .
}}\end{array}
$$ 
Therefore $F$ is an open map and hence a homeomorphism.
\end{proof}

Let $\Gamma$ be a topological group and $A$ be a CW-complex. 
A \cgagrou\ $\cali$ such that $\cali_a$ is a compact Lie group for
all $a\in A$ is called {\it proper}. When $\Gamma$ is itself a Lie group,
this is equivalent to saying that the $\Gamma$-action on 
the corresponding \spl $\Gamma$-CW-complex with isotropy groupoid 
$\cali$ is proper, see \cite[Theorem~1.23]{Lu}.
We have a classification theorem for equivariant bundles
over split $\Gamma$-spaces with proper isotropy groupoids
in \thref{Th-classi-CW}. 
First, we give a version of \lemref{rloma}. 
\begin{Lemma}\label{rloma-cell}
Let $\Gamma$ be a topological group, and $A$ be a CW-complex. 
Let $\cali$ be a \pgagrou. 
Then, any continuous representation of $\cali$ to a compact Lie group $G$ is \rloma. 
\end{Lemma}
\begin{proof}
Let $\alpha\:\cali\to G$ be a continuous representation.
Let $a\in A$. We shall construct a pair $(U,g)$, where $U$ is an open 
set of $A$, such that $\cali_u$
is a subgroup of $\cali_a$ for each $u\in U$, and $g\:U\to G$ is a continuous map
satisfying  $\alpha_u(\gamma)=g(u)\alpha_a(\gamma)g(u)^\mun$ for all
$u\in U$ and all $\gamma\in\cali_u$. Call the pair $(U,g)$ an {\em $a$-straightening}
of $\alpha$ in $A$. 
The final open set $U$ will be a neighbourhood of $a$, but
the definition of an $a$-straightening does not use that 
$a\in U$, just that the element $g(a)\in G$ is defined.
An $a$-straightening is equivalent to the data of a sequence
$(U^{(d)},g_{d})$ of $a$-straightenings of $\alpha$ in $A^{(d)}$,
such that $U^{(d+1)}\cap A^{(d)}=U^{(d)}$ and $g_{d+1}|U^{(d)}=g_{d}$.

We construct $(U^{(d)},g_{d})$ by induction on $d$, setting
$U_d=\emptyset$ if $d<d(e(a))$. If $d(e(a))=0$, we set
$U^{(0)}=\{a\}$ and $g_0(a)=1$. If $d(e(a))>0$, there exists a neighbourhood 
$U^{d(e(a))}$ of $a$ in $e(a)$ with a pointed homeomorphism 
$(U^{(d(e(a)))},a)\xrightarrow{\approx} ([-1,1]^{d(e(a))},0)$. The existence of 
$g_{d(e(a))}$ is guaranteed by \lemref{relconju}. Suppose that an $a$-straightening
$(U^{(d)},g_{d})$ of $\alpha$ in $A^{(d)}$ is constructed, with $d\geq d(e(a))$
and $a\in U_d$. For $e\in\Lambda_{d+1}(A)$,
let $\sigma_e\:(\bbd^{d+1}_e,\bbs^{d}_e)\to (A^{(d+1)},A^{(d)})$ be a characteristic map for 
the cell $e$.
Let $V_e$ be the open set of $\bbs^{d}_e$ defined by $V_e=\sigma_e^\mun(U^{(d)})$. 
Let $W_e$ be the open set of $\bbd^{d+1}_e$ defined by 
$W_e=\{tx\mid x\in V_e \hbox{ and } t\in(0,1]\}$. 
The correspondence $u\mapsto \alpha_{\sigma_e(u)}\in\hom(\cali(e),G)$ is a continuous 
representation $\alpha_e$ of the $(\Gamma,\bbd^{n+1}_e)$-groupoid $\cali(e)\times\bbd^{n+1}_e$.  
The pair $(V_e,g_d\pcirc\sigma_e)$ is a $a$-straightening of $\alpha_e$ 
in $\bbs^{n}_e$. As $W_e$ is homeomorphic to $V_e\times [0,1]$, this $a$-straightening
extends to a $a$-straightening $(W_e,g_e)$ of $\alpha_e$ in $\bbd^{n+1}_e$. 
The family $g_e$ defines a map $g_{d+1}\:U^{(n+1)}\to G$, where 
$U^{(n+1)}=\bigcup_{e\in\Lambda_{n+1}(A)}W_e$, giving rise
to the $a$-straightening of $\alpha$ in $A^{(n+1)}$. 
\end{proof}

By \lemref{rloma-cell}, the isotropy representation 
$\isor\:\nsbungt(X)\to\repi$ is defined, as in \secref{S-isor}.
The classification theorem for \spl bundles over a \spl 
$\Gamma$-CW-complex with proper isotropy groupoid takes the following
form.

\begin{Theorem}[Classification II]\label{Th-classi-CW}
Let $\Gamma$ be a Lie group, and $A$ be a CW-complex. 
Let $\cali$ be a \pgagrou.
Let $(X,\pi,\varphi)$ be a \spl $\Gamma$-CW-complex over $A$ with isotropy groupoid $\cali$.
Then, for any compact Lie group $G$, 
the isotropy representation $\isor\:\nsbungt(X)\to\repi$ is a bijection.
Moreover, any \spl bundle over $X$ is numerable.
\end{Theorem}
\begin{proof}
The proof of \thref{Th-classi-CW} is the same as that of
\thref{Th-classi}, using 
\lemref{rloma-cell} instead of \lemref{rloma}
and \proref{recon-pro-CW} instead of \proref{recon-pro}. 
Being a CW-complex, $A$ is paracompact, so each open cover is 
numerable. The last assertion of \thref{Th-classi-CW} comes
from \remref{numerable}.
\end{proof}

\begin{Remark}
The assumption that $\Gamma$ is a Lie group is only used to ensure that
the quotient projection 
$q_a\:\Gamma \to \Gamma/\cali_a$ is a (numerable) principal bundle.
If we do not care about numerability, the existence of local cross-sections
of $q_a$ holds more generally (see \cite{Mo} and \cite{Mo2}).  
\end{Remark}

As in \proref{Th-classi-ab}, \thref{Th-classi-CW} extends to a classification of all $\Gamma$-equivariant $G$-bundles if $G$ is abelian.
More precisely:

\begin{Proposition}\label{Th-classi-ab-CW} 
Let $\cali$ be a \pgagrou\ for a Lie group $\Gamma$.
Let $(X,\pi,\varphi)$ be a \spl $\Gamma$-CW-complex over $A$ with isotropy groupoid $\cali$.
Then, for any compact abelian Lie group $G$, 
one has an isomorphism of abelian groups
$$
(\isor,\varphi^*)\: \nbungt(X) \xrightarrow{\approx}
\repi\times\bung(A) \, .  
$$
Moreover, any principal $\Gamma$-equivariant $G$-bundle over $X$ is numerable.
\end{Proposition}
 
\begin{proof}
The proof \proref{Th-classi-ab-CW} is the same as that of
\proref{Th-classi-ab}, using \thref{Th-classi-CW}
instead of \thref{Th-classi}. For the numerability, observe 
that the inverse of the bijection
$(\isor,\varphi^*)$ is given by $(\isor,\varphi^*)^\mun(\xi,\eta)=
\xi\otimes\pi^*\eta$. By \thref{Th-classi}, $\xi$ is numerable. Since
$A$ is a CW-complex, $\eta$ is numerable and thus 
$(\isor,\varphi^*)^\mun(\xi,\eta)$ is numerable. Hence, any 
$\Gamma$-equivariant principal $G$-bundle over $X$ is numerable.
\end{proof}

\subsection{Examples}\label{Spltsp-ex}

\begin{ccote}\label{ex-toric}
Generalised toric manifolds of real dimension $2m$,
in the sense of \cite{DJ}, are \spl $\bbt$-spaces
where $\bbt$ is an $m$-dimensional torus.
The orbit space $A$ is a simple polytope and the section $\varphi$  is given in \cite[Lemma 1.4]{DJ}.
This includes symplectic toric manifolds, see, e.g.~\cite{Gu},
where $\pi\:X\to A$ is the moment map and $A\subset {\rm Lie}(\bbt)^*$ 
the moment polytope. Our reconstruction proposition~\ref{recon-pro}
is the topological content of Delzant's theorem \cite[Theorem~1.8]{Gu},
or \cite[Proposition~1.7]{DJ}.
\end{ccote}

\begin{ccote}\label{brihae}
When $\Gamma$ is discrete, the ``strata preserving actions with strict fundamental
domain'' of \cite[Chapter~II.12]{BH} are generalizations of \spl $\Gamma$-spaces
with a cellular isotropy groupoid. Several examples are given in 
\cite[Chapter~II.12.9]{BH}.
\end{ccote}

\rm
Several of the examples below involve \cgagrou s where $A=\Delta^{m}$ 
is the standard $m$-simplex in $\bbr^{m+1}$: 
$$
\Delta^{m}=\{(t_0,\dots,t_m)\in\bbr^{m+1}\mid t_i\geq 0 \hbox{ and }
\sum_{i=0}^mt_i=1\}\ .
$$
We use the standard simplicial structure on $\Delta^{m}$,
with $\Lambda_k(\Delta^{m})$ being the set of all subsets of
$\{0,\dots,m\}$ containing $k+1$ elements.
When $m=1,2$, we use special notations illustrated by the 
following pictures.

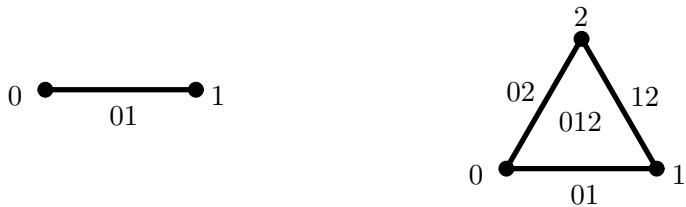
\sk{12}\small
\begin{minipage}{6cm}
\begin{center}
\setlength{\unitlength}{.03mm}
\begin{pspicture}(0,-2.1)(0,0)
\psline[linewidth=2pt](-1,0)(1,0)
\put(-1,0){\circle*{0.2}}
\put(1,0){\circle*{0.2}}
\put(-1.5,-0.2){$0$}
\put(1.2,-0.2){$1$}
\put(-.15,-0.47){$01$}
\end{pspicture}
\end{center}
\end{minipage}
\begin{minipage}{6cm}
\begin{center}
\begin{pspicture}(0,0)(0,0)
\psline[linewidth=2pt](-1,0)(1,0)(0,1.73)(-1,0)
\put(-1,0){\circle*{0.2}}
\put(1,0){\circle*{0.2}}
\put(0,1.73){\circle*{0.2}}
\put(-1.5,-0.2){$0$}
\put(-.1,1.9){$2$}
\put(1.2,-0.2){$1$}
\put(-1.0,0.9){$02$}
\put(0.65,0.83){$12$}
\put(-.15,-0.47){$01$}
\put(-0.3,0.5){$012$}
\end{pspicture}
\end{center}
\end{minipage}
\sk{-20}

\begin{ccote}\label{sphere}
Let $\bbt=(S^1)^{m+1}$. We define a \ctgrou\ $\cali$ with $A=\Delta^m$ by
$$
\cali(e)=\{(\gamma_0,\dots,\gamma_m)\mid \gamma_i=1 \hbox{ if } i\in e\} \, ..
$$
A model $(X,\pi,\varphi)$ for the \spl $\bbt$-space with isotropy groupoid
$\cali$ is given by $X=S^{2m+1}\subset\bbc^{m+1}$ with the $\bbt$-action
$(\gamma_0,\dots,\gamma_m)\cdot (z_0,\dots,z_m)=(\gamma_0\,z_0,\dots,\gamma_m\,z_m)$.
The map $\pi$ and $\varphi$ may be chosen as
\begin{equation}\label{piphi}
\begin{array}{rcl}
\pi(z_0,\dots,z_m)&=&(|z_0|^2,\dots,|z_m|^2)\
\\[1mm]
\varphi(t_0,\dots,t_m)&=&(\sqrt{t_0},\dots,\sqrt{t_m})\ .
\end{array}
\end{equation}

More generally, let $\bbt$ be any torus and let $\chi_0,\dots\chi_m\in\hom(\bbt,S^1)$.
Define a \ctgrou\ $\cali$ with $A=\Delta^m$ by $\cali(e)=\bigcap_{j\in e} \ker\chi_j$.
A model for the \spl $\bbt$-space with
isotropy groupoid $\cali$ is again given by $(S^{2m-+},\pi,\varphi)$,
where $\pi$ and $\varphi$ are defined by Equations~\eqref{piphi} and
where the $\bbt$-action on $S^{2m+1}$ is
$$
\gamma\cdot (z_0,\dots,z_m)=
(\chi_0(\gamma)z_1,\dots,\chi_m(\gamma)z_m) \ .
$$
\end{ccote}

\begin{ccote}\label{exCP}
Let $\cali$ be a \cgagrou\
and let $\Gamma_0$ be a closed subgroup of
$\Gamma$. A \cgagrou\ $\cali_0$ is then defined on $A$ by 
$\cali_0(e)$ be the
subgroup generated by $\cali(e)\cup\Gamma_0$. If $(X,\pi,\varphi)$ is the
\spl $\Gamma$-space over $A$ with isotropy groupoid $\cali$, then
that with isotropy groupoid $\cali(\Gamma_0)$ is
$(\Gamma_0\backslash X,\pi^0,\varphi^0)$, where $\pi^0$ is
induced by $\pi$ and $\varphi^0$ is $\varphi$ composed
with the projection $X\to\Gamma_0\backslash X$. For instance,
if we take $\Gamma_0$ to be the diagonal $S^1$ in \exref{sphere},
we get a \spl $(S^1)^{m+1}$-structures on the complex projective space
$\bbc P^{m}$.
\end{ccote}

\begin{ccote}\label{sonsph}
Let $\Gamma=SO(n+1)$. We see $SO(n)$ as the subgroup of $\Gamma$ leaving the last coordinate
fixed. Consider the \cgagrou\ with $A=[-1,1]$, defined by $\cali_{\pm 1}=\Gamma$
and $\cali_{(-1,1)}=SO(n)$. The \spl $\Gamma$-space with
isotropy groupoid $\cali$ is
$(S^n,\pi,\varphi)$ with $\pi(x_1,\dots,x_{n+1})=x_{n+1}$ ($\varphi$ may be defined using a meridian).
The classification of \spl $\Gamma$-equivariant $G$-bundles
over $S^n$ has been studied in \cite{hhausmann2}.
\end{ccote}

\begin{ccote}\label{coh1} Let $X$ be a 
 $\Gamma$-CW-complex $X$ so that the orbit space,
with its induced $CW$-structure, is a segment (say $\Delta^1$). 
This is one type of {\it cohomogeneity one action}. There are then
subgroups $\Gamma_0,\Gamma_1,\Gamma_{01}$ of $\Gamma$ so that $X$ is $\Gamma$-equivariantly
homeomorphic to $\Gamma/\Gamma_{01}\times [0,1]$ glued to  
$\Gamma/\Gamma_{0}\times \{0\}$ and $\Gamma/\Gamma_{1}\times \{1\}$ by equivariant maps. 
Sending $t$ to $([e],t)\in\Gamma/\Gamma_{01}\times [0,1]$ produces a splitting
$\varphi$ with a cellular isotropy groupoid $\cali$, satisfying
$\cali_{01}=\Gamma_{01}$, $\cali_{0}=\Gamma_{0}$ and $\cali_{1}=\Gamma_{1}$. 
The space  $X$ is then a \spl $\Gamma$-space with isotropy groupoid $\cali$. 
If $\Gamma$ is a compact Lie group, one checks that $X$ has a natural smooth manifold structure for which the
action is smooth. For more details and references on cohomogeneity one action,
see \cite[\S\,8]{hhausmann2}, where $\Gamma$-equivariant $G$-bundles over such 
$\Gamma$-spaces are classified (they are all \spl\hskip-4pt).
\end{ccote}

\begin{ccote}\label{ex-p1chi}\rm
Let $\bbt$ be any torus and
let $\chi$ be a non-trivial element in $\hom(\bbt,S^1)$.
Define a \ctgrou\ $\cali$ with $A=\Delta^1$ by $\cali_0=\cali_1=\bbt$ and
$\cali_{01}=\ker\chi$. The \spl $\bbt$-space with isotropy groupoid $\cali$
is $(\bbc P^1,\pi,\varphi)$, where $\pi([z_0\!\:\!z_1])=(|z_0|^2,|z_1|^2)$,
$\varphi(t_0,t_1)=[\sqrt{t_0}\:\sqrt{t_1}]$ and the $\bbt$-action
is given by $\gamma\,[x_0\!\:\!x_1]=[\chi(\gamma)\,x_0\!\:\!x_1]$.
We denote this \spl $\bbt$-space by $\bbc P^1(\chi)$.
\end{ccote}

\section{Cellular representations - Computations of $\repi$}

\subsection{Cellular representations}\label{algrep}
Let $\cali$ be a \cgagrou.  
A representation $\beta\:\cali\to G$ is called {\it cellular} 
if $\beta_a=\beta_b$ when $e(a)=e(b)$. For each $e\in\Lambda(A)$,
this thus defines $\beta_e\in\hom(\cali(e),G)$, with the face
compatibility conditions $\beta_e=\beta_f|\cali(e)$ whenever $f\leq e$.
Two cellular representations $\alpha$ and $\beta$ are called {\it conjugate}
if there exists $g\in G$ such that $\beta(\gamma)=g^\mun\alpha(\gamma)g$
for all $\gamma\in\cali_a$ and all $a\in A$.
Denote by $\arep(\cali)$ the set of conjugacy classes of cellular 
representations of $\cali$ into $G$.

To a cellular representation $\alpha\:\cali\to G$ and a cell $e$ of $A$, 
one can associate its conjugacy class $[\alpha_e]\in\chom(\cali(e),G)$. This gives
rise to a map 
$$
\kappa\: \arep(\cali)\to \prod_{e\in\Lambda(A)} \chom(\cali(e),G) \ .
$$
If an element  $(b_e)$ of this product is in the image of $\kappa$, it must
satisfy the face compatibility conditions, that is 
the equation $b_e=b_f|\cali(e)$ holds in $\chom(\cali(e),G)$
whenever $f\leq e$.
We then define
$$
\carep(\cali)=\big\{(b_e)\in \prod_{e\in\Lambda(A)}\chom(\cali(e),G)
\mid b_e=b_f|\cali(e) \hbox{ if } f\leq e \big\} 
$$
and see $\kappa$ as a map $\kappa\:\arep(\cali)\to\carep(\cali)$.
 When $\cali$ is a \pgagrou,
the map $\kappa$ sits in a commutative diagram
\begin{equation}\label{ijk}
\begin{array}{c}{\xymatrix@C-3pt@M+2pt@R-4pt{%
\arep(\cali) \ar[dr]_(0.50){\kappa}
\ar[rr]^(0.50){\jmath}  &&
\repcat \ar[dl]^(0.50){\imath}
 \\
& \carep(\cali)
}}\end{array} \ .
\end{equation}
The map $\jmath$ is obvious, since a cellular representation is
a representation, which is clearly \rloma.
To define $\imath(\beta)_e$ for $e\in\Lambda(A)$, we choose
$a\in A$ with $e(a)=e$ and set $\imath(\beta)_e = [\beta_{a}]$.
Since cells are connected, $\imath$ is
well defined by \lemref{condiscret-bis}. Although none of these
maps is either surjective or injective in general, Diagram~\eqref{ijk} is the source of all our 
information about $\repi$ so far.

One useful method for computing $\arep(\cali)$ and $\carep(\cali)$
is to restrict representations of $\cali$ to skeleta of $A$.
This yields  restriction maps 
${\rm res}_k\:\arep(\cali)\to\arep(\cali^{(k)})$ and
${\rm res}_k\:\carep(\cali)\to\carep(\cali^{(k)})$. 
Recall that a CW-complex $A$ is {\it regular} if each cell $e$ admits
a characteristic map $\sigma_e\:\bbd^{d(e)}\to A$ that is an embedding,
sending $\bbs^{d(e)-1}$ onto a subcomplex of $A^{(d(e)-1)}$.
We set $\|e\|=\sigma_e(\bbd^{d(e)})$, the closure of $|e|$.
To simplify the notations, we write $\partial e$ instead of $\partial\|e\|$
for the boundary of $\|e\|$.

\begin{Proposition}\label{res01skel} 
Let $\cali$ be a \cgagrou\ for $\Gamma$ a topological group.
Assume that $A$ is a regular CW-complex. Then, for any topological group $G$,
one has
\renewcommand{\labelenumi}{(\alph{enumi})} 
\begin{enumerate}
\item
${\rm res}_0\:\arep(\cali)\to\arep(\cali^{(0)})$ and 
${\rm res}_0\:\carep(\cali)\to\carep(\cali^{(0)})$ \\ are injective.
\item
${\rm res}_1\:\arep(\cali)\to\arep(\cali^{(1)})$ and 
${\rm res}_1\:\carep(\cali)\to\carep(\cali^{(0)})$ \\ are bijective.
\end{enumerate}
\end{Proposition}

\preu Let $\alpha\in\arep(\cali)$ (the proof for $\carep(\cali)$ is the same).
As $A$ is regular, each cell of $A$ has a face which is a vertex. 
Therefore, ${\rm res}_0(\alpha)$ determines $\alpha$ which
proves (a) and the injectivity part of (b).

For the surjectivity in (b), it is enough to prove that the restriction map
$\arep(\cali^{(n)})\to\arep(\cali^{(n-1)})$ is onto when $n\geq 2$. 
Let $\beta\:\cali^{(n-1)}\to G$ be a cellular representation. We must
extend $\beta$ to $\hat\beta=\beta\cup\{\beta_e\}\in\rep(\cali)$, which may be
done for each $n$-cell independently. For each $e\in\Lambda_n(A)$, 
choose $f\in\Lambda_{n-1}(A)$ with $f\leq e$ and define 
$\beta_e=\beta_f\,|\,\cali(e)$. We must check that $\beta_e$ does not
depend on the choice of $f$. Let $f'$ be another choice. 
As $n\geq 2$, there exists a continuous path $c(t)$
in the frontier of $|e|$ joining $a\in |f|$ to $a\in |f'|$. By the face compatibility condition,
$\beta_a|\cali_{c(t)}$ is constant, thus $\beta_f\,|\,\cali(e)=\beta_{f'}\,|\,\cali(e)$. \cqfd

\begin{Corollary}\label{res01skel-cor} 
Suppose that the hypotheses of \proref{res01skel} hold true.
Let $b\in\carep(\cali)$. If 
${\rm res}_1(b)\in\kappa(\arep(\cali^{(1)}))$,
then $b\in\kappa(\arep(\cali))$. 
\end{Corollary}

\preu
This is a consequence of the commutative diagram
\begin{equation}\label{res01skel-cor-dia}
\begin{array}{c}{\xymatrix@C-3pt@M+2pt@R-4pt{%
\arep(\cali) \ar[d]^(0.50){\kappa}
\ar[r]^(0.45){{\rm res}_1}_(0.45){\approx}  &
\arep(\cali^{(1)}) \ar[d]^(0.50){\kappa}  \\
\carep(\cali) \ar[r]^(0.45){{\rm res}_1}_(0.45){\approx}  &
\carep(\cali^{(1)})
}}\end{array} ,
\end{equation}
the bijectivity of the horizontal arrows coming from 
\proref{res01skel}. \cqfd

\subsection{Case where $G$ is abelian}

\begin{Proposition}\label{ijkabelian} 
Let $\cali$ be a \pgagrou\ for a topological group $\Gamma$.
Let $G$ be a compact abelian Lie group.
Then the three maps $\imath,\jmath,\kappa$ of Diagram~\textup{\eqref{ijk}} are
bijective.
\end{Proposition}

\preu The map $\kappa$ is bijective since conjugation has
no effect if $G$ is abelian. It is then enough to prove that
$\jmath$ is surjective. Let $\beta\in\repcat$. As in the construction of $\imath$,
one shows that $\beta(\zeta)=\beta(\zeta')$ if $e(\pi_2(\zeta))=e(\pi_2(\zeta'))$, which
is equivalent to $\beta$ being in the image of $\jmath$. 
\end{proof}

If $\Gamma$ is a Lie group, \proref{ijkabelian} 
together with the classification Theorem~\ref{Th-classi-CW}
gives a bijection $\nsbungt(X)\approx\arep(\cali)$. 
Using \proref{res01skel}
and the functorial property in Diagram~\eqref{funct-contradiag}
(which holds true in the framework of \thref{Th-classi-CW}),
this also shows that, for $G$ abelian, the restriction maps
$\nsbungt(X)\to\nsbungt(X^{(0)})$ and  $\nsbungt(X)\to\nsbungt(X^{(1)})$ 
are respectively injective and bijective,
when $A$ is a regular CW-complex.

By \lemref{Kabelian}, one has 
$G\xrightarrow{\approx} G_0\times\pi_0(G)$, where $G_0$ is the identity component of 
the unit element. Therefore, $\arep(\cali)\approx {\rm Rep}^{\pi_0(G)}_{\scriptstyle\rm cell}(\cali)
\times {\rm Rep}^{G_0}_{\scriptstyle\rm cell}(\cali)$ (the same decomposition holds for
$\nsbungt(X)$ by Equation~\eqref{prod-eq}, again true in the context of 
\thref{Th-classi-CW}). The group $G_0$ is 
isomorphic to a product of circles, so ${\rm Rep}^{G_0}_{\scriptstyle\rm cell}(\cali)$
is a product of copies of $\areps(\cali)$. We shall now study the latter.

\subsection{$\areps(\cali)$ for $\cali$ a toric groupoid.}\label{CBOBSS}
Let $\bbt$ be a torus. 
A \ctgrou\ $\cali$ is called {\it $0$-toric} if $\cali_v=\bbt$ for all
$v\in\Lambda_0=\Lambda_0(A)$. It is called {\it $1$-toric} if it is 
$0$-toric and if, for each $e\in\Lambda_1=\Lambda_1(A)$, $\cali(e)$ 
is a codimension $1$ subtorus of $\bbt$. 
There is then $\chi_e\in\hom(\bbt,S^1)$
with $\ker\chi_e=\cali(e)$.
 The part of $X$ above the closure $\|e\|$ of $|e|$
is a $\bbt$-space isomorphic to
$\bbc P^1(\chi_e)$ of Example~\ref{ex-p1chi}.
The \spl\ $\bbt$-space with isotropy groupoid $\cali^{(1)}$
is then a graph of $\bbc P^1(\chi)$'s. Such a space $\bbt$-space $X$
is also called a \emph{GKM-space}, as this
property was first studied by  M.~Goresky, R.~Kottwitz and
R.~MacPherson in \cite{GKM}.

Let $\algt$ be the Lie algebra of $\bbt$ and let $\algl=\ker({\rm exp}\:\algt\to\bbt)$.
Let $\algl^*=\{w\in\algt^*\mid w(\algl)\subset\bbz\}$ (the dual lattice).
Consider $S^1$ as $\bbr/\bbz$. The correspondence which assigns to 
$\alpha\in\hom(\bbt,S^1)$ its differential at the unit element of $\bbt$
(the {\it weight} of $\alpha$)
produces an isomorphism between $\hom(\bbt,S^1)$ and the additive group $\algl^*$. We
shall thus identify $\hom(\bbt,S^1)$ with $\algl^*$.

Let $\cali$ be a $1$-toric \ctgrou\ with $A$ a regular complex. 
By \proref{res01skel}, $\areps(\cali)$ 
injects into  $\areps(\cali^{(0)})$ which is the
direct product of character groups
\begin{equation}\label{CBOBSS-eq1}
\areps(\cali) \subset \prod_{v\in\Lambda_0} \hom(\bbt,S^1) = \prod_{v\in\Lambda_0} \algl^* \, .
\end{equation}
Let us orient each edge $e$; this determines an ordering
$\partial_-e,\partial_+e$ of the two vertices of $e$.
The character $\chi_e$ will also be seen in $\algl^*$.
A family $(a_v)_{v\in\Lambda_0}$ is said to satisfy the {\it GKM-conditions}
if, for each $e\in\Lambda_1$,
the difference $a_{\scr\partial_+e}-a_{\scr \partial_-e}$ is a multiple of $\chi_e$.
These conditions, considered in \cite{GKM}, are also discussed in
\proref{P-eqcph-imageS1} and \remref{eqcoh-rem}.

\begin{Proposition}\label{imageS1}
Let $\cali$ be a $1$-toric \ctgrou.
The image of $\areps(\cali)$ in $\prod_{v\in\Lambda_0} \algl^*$
is the set of families $(a_v)_{v\in\Lambda_0}$
satisfying the GKM-condition.
\end{Proposition}

\preu
By \proref{res01skel}, it is enough to show that
this condition characterises the image of $\rep(\cali^{(1)})$
in $\rep(\cali^{(0)})$.
Denote by $\alpha_v\in\hom(\bbt,S^1)$ the element with weight $a_v\in\algl^*$.
The three following conditions, for $e\in\Lambda_1$ are equivalent:
\renewcommand{\labelenumi}{(\alph{enumi})}
\begin{enumerate}
\item the difference $a_{\scr\partial_+e}-a_{\scr \partial_-e}$ is a multiple of $\chi_e$.
\item $\cali(e)\subset\ker \alpha_{\scr\partial_+e}\alpha_{\scr\partial_-e}^\mun$.
\item $\alpha_{\scr\partial_+e}\mid \cali(e) = \alpha_{\scr\partial_-e}\mid \cali(e)$.
\end{enumerate}
The equivalence between (a) and (b) comes from $\cali(e)$ being of
codimension 1 in $\bbt$. This proves \proref{imageS1}.
\cqfd

\begin{Example}\label{exhirz}
Let $X$ be a symplectic toric manifold of dimension $2n$.
It is a \spl $\bbt^n$-space, with $\pi\:X\to A\subset \algt^*$
being the moment map, and the isotropy groupoid $\cali$ is $1$-toric.
The moment polytope $A$ is a $n$ dimensional convex polytope of $\algt^*$.
It is known that each edge $e$ of $A$ is parallel to $\chi_e$
(see, e.g. \cite[\S~4.2.4]{Au}). By \proref{imageS1}, $\areps(\cali)$
may be visualised as the set of affine maps $\alpha\:A\to\algt^*$
such that $\alpha(\lambda_0(A))\subset\algl^*$ and
$\alpha(|e|)$ parallel to $\chi_e$ for each $e\in\lambda_1(A)$.

The left figure below shows a $2$-dimensional moment polytope for a toric manifold,
a Hirzebruch surface diffeomorphic to $\bbc P^2\sharp \overline{\bbc P^2}$.
The torus $\bbt$ is $S^1\times S^1$, $\cali_{12}=\cali_{34}=\{1\}\times S^1$,
$\cali_{14}= S^1\times\{1\}$, $\cali_{23}$ is the diagonal subgroup and
the isotropy group for the $2$-cell is trivial.
The right figure visualises two elements of $\repi$.

\hskip -20mm
\begin{minipage}{70mm}
\sk{70}
\setlength{\unitlength}{.05mm}
\begin{pspicture}(0,0)(0,0)
\psline(4.5,4.5)(6.5,4.5)(5.5,5.5)(4.5,5.5)(4.5,4.5)
\put(4.1,4.1){$1$}\put(6.6,4){$2$}\put(5.5,5.68){$3$}\put(4.1,5.65){$4$}
\put(4.8,4.8){$A$}
\end{pspicture}
\end{minipage}\hskip 40mm
\begin{minipage}{70mm}
\sk{-20}
\setlength{\unitlength}{.02mm}
\scalebox{0.7}{
\begin{pspicture}(0,0)(0,0)
\grille
\put(-0.5,-0.5){$0$}
\put(2,-3){{\pscircle*[linecolor=black]{0.15}}}
\put(4,-3){{\pscircle*[linecolor=black]{0.15}}}
\put(0,2){{\pscircle*[linecolor=black]{0.15}}}
\put(2,2){{\pscircle*[linecolor=black]{0.15}}}
\psline[linestyle=dashed](2,-3)(4,-3)(0,2)(2,2)(2,-3)
\put(1.5,-3.5){$a(v_1)$}\put(4.2,-3.5){$a(v_2)$}\put(-0.8,2.1){$a(v_3)$}
\put(2.2,2.2){$a(v_4)$}
\put(-1,-2){{\pscircle*[linecolor=black]{0.15}}}
\put(-2,-2){{\pscircle*[linecolor=black]{0.15}}}
\put(-3,-1){{\pscircle*[linecolor=black]{0.15}}}
\put(-1,-1){{\pscircle*[linecolor=black]{0.15}}}
\psline[linestyle=dashed](-1,-2)(-2,-2)(-3,-1)(-1,-1)(-1,-2)
\end{pspicture}}
\end{minipage}
\end{Example}

\sk{-14}
Let $\alpha\in\areps(\cali)$.
Let $X$ be the \spl $\bbt$-space with isotropy groupoid $\cali$.
Let $\eta$ be a \spl $\bbt$-equivariant $S^1$-principal bundle over $X$,
with isotropy representation $\alpha$. Let $e\in\Lambda_1$.
In \proref{imageS1},
the integer $k_e\in\bbz$ such that
$a_{\scr\partial_+e}-a_{\scr\partial_-e} = k_e\chi_e$
is related to the Euler number
of $\eta$ restricted to $X_e$, the part of $X$ above
the closure $\|e\|$ of $|e|$,
which is homeomorphic to $\bbc P^1$. Choose a generator $[X_e]$
of $H_2(X_e;\bbz)$. Let $\varepsilon\in H^2(X_e;\bbz)$
be the Euler class of $\eta$ restricted to $X_e$.

\begin{Proposition}\label{euler}
$\varepsilon([X_e])=\pm k_e$.
\end{Proposition}

\preu It is enough to consider the case where $X=X_e=\bbc P^1(\chi)$
for $\chi\in\hom(\bbt,S^1)$. The quotient space $A$ is then a segment,
with two $0$-cells $0$ and $1$ and a $1$-cell $e$ and we identify $A$ with
$[0,1]$. One has $\cali_0=\cali_1=\bbt$
and $\cali(e)=\ker\chi$.
The elements $\alpha_0,\alpha_1\in\hom(\bbt;S^1)$ have weights
$a_0,a_1\in\algl^*$.
The bundle $\eta$ may then be identified with the bundle
$E_\alpha\hfl{\pi}{} Y_{\cali}$ of the proof of \thref{Th-classi}.

Let $U_0=A-\{1\}$ and
$U_1=A-\{0\}$ and call $\calw_0$ and $\calw_1$ the open sets of $X$
above $U_0$ and $U_1$. One has
local sections $\sigma_i\:\calw_i\to E_\alpha$  of $\pi$
defined by $\sigma_i([\gamma,u])=[\gamma,u,\alpha_i(\gamma)]$.
Let $s\in\hom(S^1,\bbt)$ such that $\chi\pcirc s\:S^1\to\bbt/\cali_{e}$
is surjective.  Define $\hat s\:S^1\to Y_{\cali}$ by
$\hat s(\delta)=[s(\delta),1/2]$. One has
$$
\sigma_1(\hat s(\delta))=\sigma_0(\hat s(\delta))\,\alpha_0(s(\delta))^\mun
\alpha_1(s(\delta)) = \sigma_0(\hat s(\delta))\cdot\chi(s(\delta))^{\pm k_e} \, .
$$
By the classification of $S^1$-principal bundles over a $2$-sphere,
this proves Proposition~\ref{euler}.  \cqfd

\subsection{Smooth circle bundles}
Let $\cali$ be a \ctgrou\ with $A$ a regular complex.
Let $(X,\pi,\varphi)$ be a \spl $\bbt$-space with isotropy groupoid $\cali$.
Suppose that $X$ is (closed) smooth manifold and that the $\bbt$-action
is smooth. In this subsection, we relate the isotropy representation
$\isor\:\nsbungt(X)\to\repi\approx \prod_{v\in\Lambda_0} \algl^*$
with some ``moment map'' $\Phi\:X\to\algt^*$.
The material of this section is inspired by \cite{KT}.

Let $\eta = (E\hfl{p}{} X)$ be a smooth $\bbt$-equivariant \spl principal
$S^1$-bundle over $X$. Choose $\theta\in\Omega^1(E)$ be an $\bbt$-invariant
connection of the bundle $\eta$ (we see $S^1=\bbr/\bbz$, so $\lie(S^1)=\bbr$).
This gives rise to a ``moment map'' $\Phi\:X\to\algt^*$ determined as follows.
For $\xi\in\algt$, denote by $\xi_E$ the vector field on $E$ induced  by the
action of $\bbt$. The map $\Phi$ is defined by the equation
$$
\llangle{\Phi(x)}{\xi}=\theta(\xi_E(y)) \ ,
$$
for any $x\in X$ and $z\in p^\mun(x)$.
As $\theta$ is $\bbt$-invariant, the map $\Phi$ descends to a continuous map
$\bar\Phi\:A\to\algt^*$.

Let $\alpha\in\reps(\cali)$ be the isotropy representation of $\eta$. 
For each $v\in\Lambda_0(A)$, the homomorphism $\alpha_v\in\hom(\bbt,S^1)$
is determined by its weight $a_v\in\algl^*$.

\begin{Proposition}\label{mommap}
Suppose that $\cali$ is $0$-toric. Then,
for each $v\in\Lambda_0(A)$, one has $\bar\Phi(v)=a_v$.
\end{Proposition}

\preu 
Let $\xi\in\algt$. Let $v\in\Lambda_0(A)$ and $z\in E$ with $p(z)=\varphi(v)$. 
As $\varphi(v)$ is a fixed point, the vector $\xi_E(z)$ is tangent to the $S^1$-orbit
$z\cdot S^1$.
If we identify the latter with $S^1$, then $\theta(\xi_E(z))$ is the derivative of
$\alpha_v$, that is $a_v$. \cqfd 

\begin{Remark}
The figure of \exref{exhirz} suggests a possible relationship with the ``twisted polytopes''
of \cite{KT} which remains to be investigated. 
\end{Remark}

\subsection{Case where $A$ is a graph}\label{NTgraphs}

In this section, we shall determine $\repi$ for a \cgagrou\ $\cali$
when $A$ is a graph, generalising the case treated in~\cite{hhausmann2} where $A$ is a segment.
One may suppose that the graph $A$ is regular. 
Indeed, the subdivision of an edge $e$, by adding a vertex $\hat e\in|e|$ and setting
$\cali_{\hat e}=\cali(e)$, changes neither $\repi$ nor $\carep(\cali)$. 
Observe also that if $G$ is connected,  any $G$-principal bundle over $A$
is trivial, so for a \spl\ $\Gamma$-space $X$ over $A$ one has
$\nsbungt(X)=\nbungt(X)$.
We start with some preliminary material. 

\begin{Lemma}\label{nksurjtree} 
Let $\Gamma$ and $G$ be topological groups.
Let $\cali$ be a \cgagrou, where $A$ is a tree.
Then $\kappa\:\arep(\cali)\to\carep(\cali)$ is surjective.
\end{Lemma}

\preu The lemma is true for $A=\emptyset$ since then, both $\arep(\cali)$
and $\carep(\cali)$ are empty. Otherwise, let $b\in\carep(\cali)$ and
let $v$ be a vertex of $A$. Chose $\beta_v\in\hom(\cali_v,G)$ representing
$b_v$. For an edge $e$ between $v$ and $v'$, define 
$\beta_e\in\hom(\cali(e),G)$ by $\beta_e=\beta_v|\cali(e)$.
Since $b\in\carep(\cali)$, one can choose $\beta_{v'}\in\hom(\cali_{v'},G)$
which represents $b_{v'}$ such that $\beta_{v'}|\cali(e)=\beta_e$. 
This constructs a cellular representation $\beta^1$ over
the tree $A(v,1)$ of points of distance $\leq 1$ from $v$ (for the
distance where each edge has length $1$). The same methods will propagate 
$\beta^1$ to $\beta^2$, over $A(v,2)$ and then to $A(v,n)$ for all $n$.
This defines $\beta\in\arep(\cali)$ with $\kappa(\beta)=b$.
\cqfd

\begin{Lemma}\label{nksurjgraph} 
Let $\cali$ be a \cgagrou, where $\Gamma$ is a topological group 
and $A$ is a graph.
Let $G$ be a path-connected topological group.
Then $\imath\:\rep(\cali)\to\carep(\cali)$ is surjective.
\end{Lemma}

\preu 
We may suppose that $A$ is connected: otherwise, both
$\rep(\cali)$ and $\carep(\cali)$ simply decompose into 
disjoint unions over components of $A$. 
Let $A_0$ be a maximal tree of $A$ and let $\cali_0$ be the
restriction of $\cali$ over $A_0$. Let $b\in\carep(\cali)$.
By \lemref{nksurjtree}, there exists a cellular representation
$\beta\:\cali_0\to G$ such that $\kappa(\beta)=b_{|\cali_0}$.
We want to extend $\beta$ to $\hat\beta\:\cali\to G$. This can be done
by defining $\hat\beta$
over $\|e\|$ for each edge $e$ of $A\setminus A_0$. Let
$v,v'\in\Lambda_0(A_0)$ be the vertices of $e$. 
As $b\in\carep(\cali)$, there is $g\in G$ with 
$g^\mun\beta_{v}(\gamma)g=\beta_{v'}(\gamma)$ for all $\gamma\in\cali(e)$.
Since $G$ is path-connected, there exists a continuous map
$a\mapsto g_a$, from $\|e\|$ to $G$ with $g_v=1$ and $g_{v'}=g$.
For $a\in\|e\|$, we then define $\hat\beta_a\:\cali(e)\to G$ by 
$\hat\beta_a(\gamma)=g(a)^\mun\beta_{v}(\gamma)g(a)$. \cqfd 

We now introduce some material in order to describe the preimage $\imath^\mun(\alpha)$
of $\alpha\in\carep(\cali)$. 
Let $K$ be a topological group and let
$\tilde\alpha\in\hom(K,G)$. Define $\calc(\tilde\alpha)$ to be 
the centraliser of $\tilde\alpha(K)$ in $G$. Let $\tilde\alpha'\in\hom(K,G)$ be
such that $[\tilde\alpha]=[\tilde\alpha']$ in $\chom(K,G)$.
Choose $b\in G$ such that $\tilde\alpha'(\gamma)=b\tilde\alpha(\gamma)b^\mun$.
Sending $z\in\calc(\tilde\alpha)$ to $bzb^\mun$ produces
a continuous isomorphism 
$r_{\tilde\alpha',\tilde\alpha}\:\calc(\tilde\alpha)\to\calc(\tilde\alpha')$
which does not depend on the choice of $b$. Moreover, one has
$r_{\tilde\alpha'',\tilde\alpha}\pcirc r_{\tilde\alpha'',\tilde\alpha'}=
r_{\tilde\alpha'',\tilde\alpha}$. Therefore, a topological group
$\calc(\alpha)$ is defined for $\alpha\in\chom(K,G)$: take the disjoint union
of $\calc(\tilde\alpha)$ for all representatives $\tilde\alpha$ of $\alpha$
and identify $z\in\calc(\tilde\alpha)$ with 
$r_{\tilde\alpha',\tilde\alpha}(z)\in\calc(\tilde\alpha')$. 
If $K'$ is a subgroup of $K$, one checks that $\calc(\alpha)$ is a subgroup of 
$\calc(\alpha_{|K'})$

Let $\cali$ be a \cgagrou\ and $\alpha\in\carep(\cali)$.
Let $\dot A$ be the first barycentric subdivision of $A$. 
We assume that $A$ is regular, so $\Lambda_1(\dot A)$ is the set of pairs 
$(v,e)\in\Lambda_0(A)\times\Lambda_1(A)$
with $v<e$; the edge corresponding to $(v,e)$ joins $v$ to the
barycentre $\hat e$ of $e$.
Form the group $X(\alpha)=\prod_{(v,e)\in\Lambda_1(\dot A)}\pi_0(\calc(\alpha_e))$.
Let $J^0\:\prod_{v\in\Lambda_0(A)}\pi_0(\calc(\alpha_v))\to X(\alpha)$ be the homomorphism 
sending $(x_v)$ to $(z_{(w,e)})$ with $z_{(w,e)}=j_{w,e}(x_w)$, where 
$j_{w,e}\:\pi_0(\calc(\alpha_w))\to\pi_0(\calc(\alpha_e))$ is the homomorphism induced
by the inclusion. Consider also the homomorphism \\ 
$J^1\:\prod_{e\in\Lambda_1(A)}\pi_0(\calc(\alpha_e))\to X(\alpha)$
sending $(y_e)$ to $(z_{(w,f)})$, where $z_{(w,f)}=y_f$. Set $Y^0(\alpha)$ 
and $Y^1(\alpha)$ to be
the images of $J^0$ and $J^1$ and consider the double coset family
 $Z(\alpha)=Y^0(\alpha)\backslash X(\alpha)/Y^1(\alpha)$.

\begin{Theorem}\label{NTgraph} 
Let $\cali$ be a \pgagrou, with $\Gamma$ a topological group and $A$ a graph. 
Let $G$ be a compact connected Lie group. Then, $\imath\:\repi\to\carep(\cali)$
is surjective and 
the preimage $\imath^\mun(\alpha)$ of $\alpha\in\carep(\cali)$ is in
bijection with $Z(\alpha)$. 
\end{Theorem}

Before proving \thref{NTgraph}, we state some of its corollaries, in which
we assume the hypotheses of \thref{NTgraph} and mention only the additional
hypotheses.

\begin{Corollary}\label{NTgraph-cor1}
Suppose that $A$ is a finite graph. Then,
the preimages of $\imath\:\repi\to\carep(\cali)$
are finite.
\end{Corollary}

\preu As $\calc(\alpha_e)$ is a closed subgroup in $G$, $\pi_0(\calc(\alpha_e))$ is finite
for each edge $e$ of $A$. Therefore $Z(\alpha)$ is finite. \cqfd 

The next corollary corresponds to \cite[Theorem B and 8.12]{hhausmann2}.

\begin{Corollary}\label{NTgraph-cor2}
Let $\cali$ be a proper $(\Gamma,\Delta^1)$-groupoid. Then, 
the preimage $\imath^\mun(\alpha)$ of $\alpha\in\carep(\cali)$ is in
bijection with the set of double cosets
$\pi_0(\calc(\alpha_0))\backslash\pi_0(\calc(\alpha_{01}))/\pi_0(\calc(\alpha_1))$.
\end{Corollary}

\preu The group $X(\alpha)$ is isomorphic to 
$\pi_0(\calc(\alpha_{01}))\times\pi_0(\calc(\alpha_{01}))$ with 
$Y^1(\alpha)\approx\pi_0(\calc(\alpha_{01}))$ being the diagonal subgroup.
The group $Y^0(\alpha)$ is $\pi_0(\calc(\alpha_0))\times\pi_0(\calc(\alpha_1))$.
Therefore, the map $X(\alpha)\to\pi_0(\calc(\alpha_{01}))$ given by 
$(z_0,z_1)\mapsto z_0 z_1^\mun$ descends to a bijection from
$Z(\alpha)$ to 
$\pi_0(\calc(\alpha_0))\backslash\pi_0(\calc(\alpha_{01}))/\pi_0(\calc(\alpha_1))$
(see \cite[Section 8]{hhausmann2}). \cqfd

\begin{Corollary}\label{NTgraph-cor3}
Suppose that $\cali(e)$ is a torus for all edges $e$ of $A$. Then
$\imath\:\repi\to\carep(\cali)$ is a bijection.
\end{Corollary}

\preu Let $\alpha\in\carep(\cali)$ and let  $\tilde\alpha\:\cali\to G$ be a continuous representation with $\imath(\tilde\alpha)=\alpha$
Our hypotheses imply that $\tilde\alpha_a(\cali_e)$ is a torus
for all $e\in\Lambda_1(A)$ and all $a\in|e|$. 
As $G$ is connected,
the group $\calc(\alpha_e)$ is then connected (see, e.g. \cite[Theorem~3.3.1]{DK}). 
Therefore, $Z(\alpha)$ reduces to a single element. \cqfd

\sk{1}
\ppreu{\thref{NTgraph}} The surjectivity of $\imath$ is established in \lemref{nksurjgraph}. 
Let $\alpha\in\carep(\cali)$. 
The strategy is to construct a transitive action of $X(\alpha)$ on $\imath^\mun(\alpha)$
and study the stabilisers.

Let $\tilde\alpha^0\:\cali^{(0)}\to G$ be a representative of $\alpha^{(0)}$.
Let $\tilrep(\cali,\tilde\alpha^0)$ be the set of continuous
representations from $\cali$ to $G$ which restrict to $\tilde\alpha^0$
on $\alpha^{(0)}$. As $G$ is connected, any map from $A^{(0)}$ to $G$ 
extends to $A$, which implies that each class in $\imath^\mun(\alpha)$
has a representative in $\tilrep(\cali,\tilde\alpha^0)$.
Also, if $\tilde\alpha\in\tilrep(\cali,\tilde\alpha^0)$, then $\imath(\tilde\alpha)=\alpha$
by \proref{res01skel}. Thus, the map $\tilde\alpha\mapsto[\tilde\alpha]\in\repi$
produces a surjection $\tilrep(\cali,\tilde\alpha^0)\onto\imath^\mun(\alpha)$.

Form the group $\tilde X(\tilde\alpha^0)=\prod_{(v,e)\in\Lambda_1(\dot A)}\calc(\tilde\alpha^0(\cali(e))$.
Let $z=(z_{(v,e)})\in\tilde X(\tilde\alpha^0)$ and $\tilde\alpha\in\tilrep(\cali,\tilde\alpha^0)$.
For each edge $e$ of $A$ with $\partial e = \{v,v'\}$, choose, using that $G$ is connected, 
a continuous map $g_e\:\|e\|\to G$ such that $g_e(v)=z_{(v,e)}$ and $g_e(v')=z_{(v',e)}$.
We call $\{g_e\}$ a {\it connecting family} for $z$.
Define $z\cdot_{\{g_e\}}\tilde\alpha\in\tilrep(\cali,\tilde\alpha^0)$ by
\begin{equation}\label{NTgraph-eq10}
z\cdot_{\{g_e\}}\tilde\alpha(\gamma)=
\left\{ \begin{array}{lll}
g_e(a)\tilde\alpha(\gamma)g_e(a)^\mun & 
\hbox{ if  $a\in|e|$ and $\gamma\in\cali_a$ } \\
\tilde\alpha_a(\gamma) & \hbox{ otherwise}\, .
\end{array}\right.
\end{equation}
For two connecting families $\{g_e\}$ and $\{\bar g_e\}$ for $z$, we check that
$$
z\cdot_{\{g_e\}}\tilde\alpha(\gamma)=
h(a)\big(z\cdot_{\{\bar g_e\}}\tilde\alpha(\gamma)\big)h(a)^\mun \ ,
$$ where
$h\:A\to G$ is the (continuous) map defined by
$h(a)=g_e(a)\bar g_e(a)^\mun$ if $a\in\|e\|$.
This thus defines $z\cdot\tilde\alpha$ in $\imath^\mun(\alpha)$ which does not depend on the
choice of the connecting family $\{g_e\}$.

Now, suppose that $\tilde\alpha,\tilde\alpha'\in\tilrep(\cali,\tilde\alpha^0)$
represent the same element in $\repi$. This means that there is a map
$h\:A\to G$ such that $\tilde\alpha'_a(\gamma)=h(a)\tilde\alpha_a(\gamma)h(a)^\mun$.
Observe then that $h(v)\in\calc(\tilde\alpha^0(\cali_v))$ for all $v\in\Lambda_0(A)$ and hence
$$
h(a)\big(z\cdot_{\{g_e\}} \tilde\alpha\big)h(a)^\mun = z\cdot_{\{h(a)g_eh(a)^\mun\}} \tilde\alpha' \, .
$$
We have thus defined an action of $\tilde X(\tilde\alpha^0)$ on $\imath^\mun(\alpha)$.
We now prove that this action is transitive. 
Let $\tilde\alpha,\tilde\alpha'\in\tilrep(\cali,\tilde\alpha^0)$. Orient each
edge $e$ of $A$, getting then $\partial e =\{\partial_-e,\partial_+e\}$. By \lemref{relconju}, there exist
$s,s'\:\|e\|\to G$ such that $\tilde\alpha(\gamma)=s(a)^\mun\tilde\alpha^0_{\partial_-e}(\gamma)s(a)$ and
$\tilde\alpha'(\gamma)=s'(a)^\mun\tilde\alpha^0_{\partial_-e}(\gamma)s'(a)$ for all $\gamma\in\cali(e)$ and all
$a\in\|e\|$. This implies $s(\partial e)$ and $s'(\partial_-e)$ are contained in $\calc(\tilde\alpha(\cali(e)))$.
Hence,  one has $\tilde\alpha'=z_{\{g_e\}}\cdot\tilde\alpha$, where 
$z_{(\partial_-e,e)}=s'(\partial_-e)^\mun s(\partial_-e)$, 
$z_{(\partial_+e,e)}=s'(\partial_+e)^\mun s(\partial_+e)$ and $g_e(a)=s'(a)^\mun s(a)$. 
Hence, the action of $\tilde X(\tilde\alpha^0)$ on $\imath^\mun(\alpha)$
is transitive.

If $z=(z_{(v,e)})$ is in the unit component of $\tilde X(\tilde\alpha^0)$, then 
$z\cdot_{\{g_e\}}\tilde\alpha=\tilde\alpha$,
if the maps $g_e$ are chosen so that $g(\hat e)=1$ and
$g_e(\|(v,e)\|)\subset\calc(\tilde\alpha^0_v(\cali(e)))$.  
This implies that
the action of $\tilde X(\tilde\alpha^0)$ on $\imath^\mun(\alpha)$ descends to an action of
the group $\prod_{(v,e)\in\Lambda_1(\dot A)}\pi_0(\calc(\tilde\alpha^0_v(\cali(e)))$ which is isomorphic to $X(\alpha)$.

Let $f$ be an edge of $A$, with $\partial f = \{v,v'\}$. 
The representation $\tilde\alpha^0\:\cali^{(0)}\to G$ can be chosen such that the restrictions to $\cali_e$
of $\tilde\alpha^0_v$ and $\tilde\alpha^0_{v'}$ coincide. 
For each $\zeta\in\calc(\tilde\alpha^0_v(\cali(f))=\calc(\tilde\alpha^0_{v'}(\cali(f))$
we can then consider the element $z(\zeta)$ of 
$\tilde X(\tilde\alpha^0)$ satisfying $z_{(v,f)}(\zeta)=z_{(v',f)}(\zeta)=\zeta$ and $z_{(w,e)}(\zeta)=1$ of $e\neq f$.
Then $z(\zeta)\cdot_{\{g_e\}}\tilde\alpha=\tilde\alpha$ if the $g_e$ are constant maps. This may be done for each
edge $f$ of $A$, showing that the group  $Y^1(\alpha)$ acts trivially on $\beta$ for all $\beta\in\imath^\mun(\alpha)$.

Let $y\in\prod_{v\in\Lambda_0(A)}(\calc(\tilde\alpha^0_v(\cali_v)))$ and  
$z\in\tilde X(\tilde\alpha^0)$. Consider the element $yz\in\tilde X(\tilde\alpha^0)$ defined by
$(yz)_{(v,e)}=y_vz_{(v,e)}$. Choose a connecting family $g_e\:\|e\|\to G$ for $z$.
For each $(v,e)\in\Lambda_1(\cdot A)$, choose $h_{(v,e)}\:\|(v,e)\|\to G$ such that 
$h_{(v,e)}(v)=y_v$ and $h_{(v,e)}(\hat e)=1$. This defines a continuous map $h\:A\to G$, by
$h(a)=h_{(v,e)}(a)$ if $a\in\|(v,e)\|$, which conjugates $(yz)\cdot_{\{hg_e\}}\tilde\alpha$
with $z\cdot_{\{g_e\}}\tilde\alpha$. This shows that $ux\cdot\beta=x\cdot\beta$ in $\imath^\mun(\alpha)$,
for all $u\in Y^0(\alpha)$, $x\in X(\alpha)$ and $\beta\in\imath^\mun(\alpha)$.

Fix $\beta\in\imath^\mun(\alpha)$, represented by $\tilde\beta\in\tilrep(\cali,\tilde\alpha^0)$. 
Consider the map $\tilde\Psi\:\tilde X(\tilde\alpha^0)\to\imath^\mun(\alpha)$ given by 
$\tilde\Psi(z)=[z\cdot\tilde\beta]$. 
By the above, we have shown that $\tilde\psi$ descends to
a surjection $\Psi\:Z(\alpha)\onto\imath^\mun(\alpha)$. It remains to show that $\Psi$ is
injective.  Let $z'\in\tilde X(\tilde\alpha^0)$ 
with $\Psi(z')=\Psi(z)$. Choose connecting families $\{g_e\}$ and $\{g'_e\}$ for $z$ and $z'$.
If $\Psi(z')=\Psi(z)$, there exists a map $h\:A\to G$ with 
$(z'\cdot_{\{g_e'\}}\tilde\beta)(\gamma)=h(a)(z\cdot_{\{g_e\}}\tilde\beta)(\gamma)h(a)^\mun$. Observe that
$h(v)\in\calc(\tilde\alpha^0(\cali_v))$ and therefore $h^{(0)}\:\cali^{(o)}\to G$ defines an element 
$y\in\prod_{v\in\Lambda_0(A)}(\calc(\tilde\alpha^0_v(\cali_v)))$ satisfying
$((yz)\cdot_{\{hg_e\}}\tilde\beta)(\gamma)=h(a)(z\cdot_{\{g_e\}}\tilde\beta)(\gamma)h(a)^\mun$.
Let $\bar z = yz$ and $\bar g_e=h g_e$. One has $[\bar z]=u[z]$ in $X(\alpha)$ with $u\in Y^0(\alpha)$.
Thus, $\bar z$ and $z$ represent the same class in $Y^0(\alpha)\backslash X(\alpha)$ and the equality
$z'\cdot_{\{g_e'\}}\tilde\beta =\bar z\cdot_{\{\bar g_e\}}\tilde\beta$ holds in $\tilrep(\cali,\tilde\alpha^0)$. 
Therefore, $\bar g_e(a)^\mun g'_e(a)\in\calc(\beta_a(\cali(e))$ for all $a\in\|e\|$. This implies
that $z'$ and $\bar z$ represent the same class in $X(\alpha)/Y^1(\alpha)$. Finally, we have shown
that $z$ and $z'$ represent the same class in $Z(\alpha)$, proving the injectivity of $\Psi$.
\cqfd

\sk{2}
We now give some examples of the use of \thref{NTgraph}.

\begin{ccote}\rm
Let $A$ be the $1$-simplex $\Delta^1$.
Let $\Gamma=SO(n)$, with $n=2k+1\geq 3$ and consider
the \cgagrou\ $\cali$ with
$\cali_0=\cali_1=\Gamma$ and $\cali_{01}=SO(n-1)$ (the \spl $\Gamma$-space $X$
with isotropy groupoid $\cali$ is $S^n$ with the $SO(n)$-action
fixing the north and the south pole). For $G=SO(n)$, $\carep(\cali)$
contains two elements, the trivial representation and the representation $\alpha$
with $\alpha_0=\alpha_1={\rm id}$. The preimage by $\imath$ of the trivial 
representation contains one element but $\imath^\mun(\alpha)$ contains 
two elements.  For details and developments, see \cite[Example~7.5]{hhausmann2}. 
\end{ccote}

\begin{ccote}\rm
 If $\pi_1(G)=\{1\}$, \thref{NTgraph} extends to a \cgagrou\ $\cali$ where
$A$ is of dimension 2, provided $\cali(e)=\{1\}$ when $e\in\Lambda_2(A)$.
Examples are given by toric manifolds of (real) dimension $4$. 
\end{ccote}

\begin{ccote}\rm
Let $\cali$ be the $S^1$-structure on the $1$-simplex $\Delta^1$
with $\cali_0=\cali_1=S^1$ and $\cali_{01}=\{1\}$. The \spl $S^1$-space with
isotropy groupoid $\cali$ is $S^2$ with $S^1$ acting by rotation
around an axis. By \thref{NTgraph},
${\rm SBun}^{G}_{S^1}(S^2)\approx \chom(\cali_0,G)\times \chom(\cali_1,G)$.
Choosing a maximal torus $\calt$ in $G$, this yields
${\rm SBun}^{G}_{S^1}(S^2)\approx \hom(\cali_0,\calt)/\calw\times \hom(\cali_1,\calt)/\calw$
where $\calw$ is the Weyl group for $\calt$.
If $G$ is of rank $k$, then $\hom(\cali_0,\calt)$ and $\hom(\cali_1,\calt)$
are both in bijection with $\bbz^k$.

Let us specialise to $G=SO(m)$ for $m\geq 3$.
A maximal torus $\calt$ of $SO(m)$ is formed by matrices containing 2-blocks
concentrated around the diagonal, so isomorphic to $SO(2)^k$,
and where $k=[m/2]$. The action of $\calw$ on $\hom(S^1,\calt)\approx\bbz^k$
can be deduced from \cite[p.~114]{Ad}.
When $m=2k+1$, the action of $\calw$ on $\bbz^k$ is generated by the
permutation of coordinates and sign changes in any of them.
A fundamental domain $\cald\subset\bbz^k$ is then
$$\cald=\{(r_1,\dots,r_k)\in\bbz^k\mid 0\leq r_1 \leq\cdots\leq r_k  \}$$
and ${\rm SBun}^{SO(2k+1)}_{S^1}(S^2)\approx \cald\times\cald$.
When $m=2k$, the sign changes must be even in number.
A fundamental domain $\cale\subset\bbz^k$ is then
$$\cale=\{(r_1,\dots,r_k)\in\bbz^k\mid
0\leq r_1 \leq\cdots\leq r_{k-1}\leq |r_k|\}$$
and ${\rm SBun}^{SO(2k)}(S^2)_{S^1}\approx \cale\times\cale$.

This example was treated in our paper \cite[Example~7.3]{hhausmann2} but the
determination of ${\rm SBun}^{SO(m)}_{S^1}(S^2)$ is wrong there
because, in the action of the Weyl group, the sign changes were
forgotten. However, the computation in \cite[Example~7.3]{hhausmann2} of the second Stiefel-Whitney number $w_2(\xi)$ for $\xi\in {\rm SBun}^{SO(m)}(S^2)_{S^1}$,
being ${\rm mod\,} 2$, is correct.
\end{ccote}

Here is an interesting consequence of the proof of \thref{NTgraph}.

\begin{Proposition}\label{jinjective} 
Let $\cali$ be a \pgagrou, with $A$ a regular CW-complex and $\Gamma$ a topological group. 
Let $G$ be a compact connected Lie group. Then $\jmath\:\arep(\cali)\to\repi$ is injective.
\end{Proposition}

\preu
As in the proof of \lemref{nksurjgraph}, one may assume that $A$ is connected.
Let $A_0$ be a maximal tree of $A$ and let $\cali_0$ be the restriction 
of $\cali$ over $A_0$.
As $A$ is connected, $A_0$ contains all the vertices of $A$ and then the restriction map
$\arep(\cali)\to\arep(\cali_0)$ is injective by \proref{res01skel}. Therefore, it is enough to prove \proref{jinjective} when $A$ is a tree. 

Let $\beta,\beta'\:\cali\to G$ be cellular representations with $\jmath(\beta)=\jmath(\beta')$. 
Let $v$ be a vertex of the tree $A$. By conjugation of $\beta$ with a constant element of $G$,
one may assume that $\beta_v=\beta'_v$. Let $e$ be an edge between $v$ and $v'$; one has
$\beta_e=\beta'_e$. Suppose that $\beta_{v'}\neq\beta'_{v'}$. Then $\beta'_v(\gamma)=z\beta(\gamma)z^\mun$ with 
$z\in\calc(\beta_{e}(\cali(e))$ and $z\notin\calc(\beta_{v'}(\cali_{v'}))$. Choose a continuous map
$g_e\:\|e\|\to G$ with $g_e(v)=1$ and $g_e(v')=z$. Let $\cali_{\|e\|}$ be the restriction of 
$\cali$ over $\|e\|$ and let $\beta''\:\cali_{\|e\|}\to G$ be the (non-cellular) representation defined by 
$\beta''(\gamma)=g_e(a)^\mun\beta'(\gamma)g_e(a)$. 
Using the notations of the proof of \thref{NTgraph}, this means that
$\beta,\beta''\in\tilrep(\cali_{\|e\|},\beta^{(0)})$
and $\beta''=y\cdot_{g_e}\beta$, where $y\in\tilde X(\beta^{(0)})$ is defined by $y_{(v,e)}=1$ and 
$y_{(v',e)}=z$. The element $y$ is non-trivial in $Z(\alpha_{\|e\|})$ which, by \thref{NTgraph},
would contradict the assumption $\jmath(\beta)=\jmath(\beta')$. Therefore, $\beta_{v'}=\beta'_{v'}$.
This argument may be done independently for all edges adjacent to $v$ and then propagated to the
whole tree $A$. \cqfd

When $A$ is a tree, the map $\jmath\:\arep(\cali)\to\repi$ is actually bijective. More precisely:

\begin{Lemma}\label{treegraph} 
Let $\cali$ be a \pgagrou, where $A$ is a graph and $\Gamma$ a topological group.
Let $A_0$ be a subtree of $A$. Let $G$ be a compact Lie group.
Then, any $\alpha\in\repi$ has a representative which is cellular over $A_0$. 
\end{Lemma}

\preu Let $v$ be a vertex of $A_0$. For an edge $e$ of $A_0$, between $v$ and $v'$, 
there exists, by \lemref{relconju}, a map $\psi_e\:\|e\|\to G$ such
that $\psi_e(a)\alpha_a(\gamma)\psi_e(a)^\mun=\alpha_v(\gamma)$ for each $a\in\|e\|$
and $\gamma\in\cali_e$. This defines a map $\psi_1\:A_0(v,1)\to G$
(notations as in the proof of \lemref{nksurjtree}). As $A_0(v,1)$ is contractible,
the homotopy extension property permits us to extend $\psi_1$ to a continuous map
$\psi_1\:A\to G$. The maps $\psi_1$ conjugates $\alpha$ to $\alpha_1$ which
is cellular over $A_0(v,1)$. The process propagates over $A_0(v,n)$ for all $n$,
giving rise to a map $\psi\:A\to G$ which conjugates $\alpha$ to a representation
which is cellular over $A_0$. \cqfd

\proref{jinjective} together with \lemref{treegraph} imply the following

\begin{Corollary}\label{jinjective-cor} 
Let $\cali$ be a \pgagrou\ with  $\Gamma$ a topological group and $A$ a tree.
Let $G$ be a compact connected Lie group. Then $\jmath\:\arep(\cali)\to\repi$ is bijective. 
\end{Corollary}
 
\begin{ccote}\label{kappanotonto}
In contrast with \thref{NTgraph}, the map $\jmath\:\arep(\cali)\to\repi$
is not surjective when the graph $A$ is not a tree. 
Using \lemref{nksurjtree}, it is enough to find an example where
$\kappa\:\arep(\cali)\to\carep(\cali)$ is not surjective.
Let $A$ be the $1$-skeleton of the $2$-simplex $\Delta^2$ with
$\cali_0=\cali_1=\cali_2=\Gamma=S^1\times S^1$, $\cali_{01}=1\times S^1$,
$\cali_{02}=S^1\times 1$ and $\cali_{12}$ is the diagonal $S^1$.
The \spl $\Gamma$-space with this isotropy groupoid is $\bbc P^2$ with
the action $(c_1,c_2)\cdot[z_0\:z_1\:z_2)]=[c_0z_0\:c_1z_1\:z_2)]$.
Take $G=SU(2)$; the diagonal torus $H$ has dimension $1$ and its Weyl group
$\calw$ acts by passing to the inverse. Then
$$\chom(\Gamma,SU(2))\approx \hom(\Gamma,H)/\calw\approx
\hat\Gamma/\{\chi\,\sim\, -\chi\}
\approx (\bbz\times\bbz)/\{(p,q)\,\sim\,-(p,q)\}.$$
We identify $\chom(\Gamma,SU(2))$
with the fundamental domain $\cald$ in $\bbz\times\bbz$:
$$\cald\:=\{(p,q)\in\bbz\times\bbz)\mid q\geq 0 \hbox{ and }
(p\geq 0 \hbox{ if } q=0)\}.$$

If $\beta\in\arep(\cali)$ is not trivial, it must be not trivial on
at least one edge-isotropy groups (say $\cali_{01}$). Then $\beta$ is conjugate
to $\beta'$ such that $\beta_{01}'(\cali_{01})\subset H$.
As $H$ is maximal abelian, $\beta'$ is then an algebraic representation of
$\cali$ in $H$. By \proref{imageS1}, one has an identification
of ${\rm Rep}^H(\cali)$ with the set of triples
$$\big((p_0,q_0),(p_1,q_1),(p_2,q_2) \big)\in (\bbz\times\bbz)^3$$
such that 
\begin{equation}\label{kappanotonto-eq1}
p_0=p_2 \ ,\ q_0=q_1 \ \hbox{ and }p_1+q_1=p_2+q_2 \ .
\end{equation}
A class in $\carep(\cali)$ is a triple
$$([p_0,q_0],[p_1,q_1],[p_2,q_2])\in\cald\times\cald\times\cald$$
such that $|p_0|=|p_2|$, $q_0=q_1$ and $|p_1+q_1|=|p_2+q_2|$.
The class $\alpha\in\carep(\cali)$ corresponding to
$([-1,2],[3,2],[1,4])$ is not in the image of $\kappa$.
Indeed, none of the 8 triples in
$(\bbz\times\bbz)^3$ above $\alpha$ satisfies 
Equations~\eqref{kappanotonto-eq1}. 
\end{ccote}

\section{Comparison with the homotopy-theoretic approach}\label{shapproach}

\subsection{Haefliger classifying spaces}\label{haefborel}

Let $(X,\pi,\varphi)$ be a \spl $\Gamma$-space over a space $A$
with isotropy groupoid $\cali$. Let $B\cali$ be the Haefliger
classifying space for $\cali$ \cite[p.~140]{Hae1}.
For a groupoid like $B\cali$ where morphisms go from one object to itself,
we check that the construction of \cite[p.~140]{Hae1} takes the following form:
set
$$
E\cali = \{(v,a)\in E\Gamma\times A\mid v\in E\cali_a\} \, ,
$$
with the induced topology, and define $B\cali$ as the  
quotient space $E\cali/\cali$. The projection $\bar\pi\:B\cali\to A$
makes $B\cali$ is a space over $A$
whose stalk over $a$ is the Milnor classifying space $B\cali_a$.
There is a section $j\:A\to B\cali$ of $\bar\pi$, sending $a\in A$ to
the class of $(v_0,a)$ where $v_0=(1e,0,\dots)\in E\Gamma$, expressed as the infinite join, 
with $e$ the unit element of $\Gamma$. 
The inclusion $\cali\subset\Gamma\times A$ is a morphism of topological
groupoids and therefore induces a continuous map $E\cali\to E\Gamma\times A$
which descends to a continuous map $B\cali\to B\Gamma\times A$.

Recall that the {\it Borel construction} associates to $X$
the space $X_\Gamma=E\Gamma\times_\Gamma X$.
The map $\pi\:X\to A$ descends to a continuous and open surjective map
$\bar\pi\:X_\Gamma\to A$, with $\bar\pi^\mun(a)=E\Gamma\times_{\Gamma}\cali_a
\approx B\cali_a$.
The composed map $E\cali\to E\Gamma\times A \hfl{{\rm id}\times \varphi}{}
E\Gamma\times X$ descends to a continuous map
$\delta\:B\cali\to X_\Gamma$ over the identity of $A$. The restriction of $\delta$
to each stalk is a weak homotopy equivalence.
It would then be interesting to figure out, for instance in the spirit of Sections~\ref{Ssplsp}
and~\ref{nspeqbd}, under which hypotheses $\delta$ is a weak homotopy equivalence.
We will restrict ourselves to \cgagrou s, where we get the following proposition.

\begin{Proposition}\label{borel-sp} 
Let $\cali$ be a \pgagrou\ for a Lie group $\Gamma$. 
Let $(X,\pi,\varphi)$ be a \spl $\Gamma$-space over $A$,
with isotropy groupoid $\ciat$. Then, the
map $\delta\:B\ciat\to X_\Gamma$ is a homotopy equivalence.
\end{Proposition}

\preu 
By \proref{recon-pro-CW}, we may suppose that 
$(X,\pi,\varphi)=(Y_{\cali},\varpi,\phi)$.
If $K$ is a subspace of $A$, we denote by $\cali(K)$
the subgroupoid of $\cali$ formed by all the stalks over $K$,
and we set $X(K)=Y_{\cali(K)}$.

Observe first that \proref{borel-sp} is true if $\cali_a$ is constant
for all $a\in A$. Indeed, one then has 
$X=\Gamma/\cali_a \times A$, so $B\cali\approx B\cali_a\times A$ and
$X_\Gamma\approx E\Gamma/\cali_a\times A$ and $\delta$ 
is a homotopy equivalence. More generally, \proref{borel-sp} 
remains true if $\cali_a$ is locally constant, meaning constant on each connected
component of $A$.

\proref{borel-sp} will be proved, by induction on $n$, for  
$X(A^{(n)})$, the
\spl $\Gamma$-space over the $n$-skeleton $A^{(n)}$ of $A$.
It is true for $n=0$ since $\cali(A^{(0)})$ is locally constant.
The induction step involves the
subcomplexes $K'=A^{(n-1)}\subset K=A^{(n)}$, 
so $K$ is obtained from $K'$ by
adjunction of $\cale=\coprod_{e\in\Lambda_n}D^n_e$, via the attaching map
$f\:\partial\cale=\coprod_{e\in\Lambda_n}S^{n-1}_e\to K'$ ($\Lambda_n=\Lambda_n(A)$). 
Then, $X(K)$ is obtained
from $X(K')$ by attaching the $\Gamma$-space  
$\tilde\cale=\coprod_{e\in\Lambda_n}(\Gamma/\cali(e)\times D^n_e)$
via the $\Gamma$-equivariant map 
$\tilde f\:\partial\tilde\cale=\coprod_{e\in\Lambda_n}(\Gamma/\cali(e)\times S^{n-1}_e)\to X(K')$.
We denote by $F\:\cale\to K$ and $\tilde F\:\tilde\cale\to X(K)$ 
the characteristic maps, extending $f$ and $\tilde f$.  
We see $\tilde\cale$ and $\partial\tilde\cale$ as \spl $\Gamma$-spaces
over $\cale$ and $\partial\cale$ respectively with locally constant
isotropy groupoids: if $x\in D^n_e$, 
then $\cali(\cale)=\cali(\partial\cale)=\cali(e)$. 
Let us consider the following diagram:
$$
\xymatrix@C-5pt@R-4pt{%
B\cali(\partial\cale) \kern 3pt  \ar[ddd]_{Bf}
\ar@{>->}[rrr] \ar@{>->}[dr]^{\simeq}_{\delta_{\partial\cale}}
&\ar@{} [dr] |{\textbf{II}}&&
B\cali(\cale) \ar[ddd]^{BF} \ar@{>->}[dl]_{\simeq}^{\delta_{\cale}}
\\
\ar@{} [dr] |{\textbf{I}}& 
(\partial\tilde\cale)_\Gamma \kern 3pt  \ar@{>->}[r]  \ar[d]_{\tilde f_\Gamma}
& \tilde\cale_\Gamma \ar[d]_{\tilde F_\Gamma}
\ar@{} [dr] |{\textbf{III}}&
\\
& X(K')_\Gamma \kern 1mm \ar@{>->}[r]\ar@{} [dr] |{\textbf{IV}}
& X(K)_\Gamma&
\\
B\cali_{K'} \kern 1mm\ar@{>->}[rrr] \ar@{>->}[ur]^{\delta_{K'}}_{\simeq}
&&&
B\cali_{K} \ar[ul]_{\delta_K}
}
$$
\noindent
The maps $\delta_{\partial\cale}$ and $\delta_{\cale}$ are homotopy
equivalences since the isotropy groupoids are locally constant.
The map $\delta_{K'}$ is a homotopy equivalence by induction hypothesis.
Restriction to any stalk shows that Diagrams I--IV are commutative.
As $\Gamma$ is a Lie group, all the spaces under consideration
have the homotopy type of CW-complexes. Therefore,
the outer and inner square diagrams
are homotopy push-out diagrams. By push-out properties,
the map $\delta_{K}$ is a homotopy equivalence. 
\cqfd

\subsection{Split bundles and classifying spaces}\label{spbunborel}

Let $\eta\:(P\xrightarrow{p} X)$ be a $\Gamma$-equivariant
principal $G$-bundle.
The Borel construction $E\Gamma\times_\Gamma P\to E\Gamma\times_\Gamma X$
yields a principal $G$-bundle $\eta_{\scr \Gamma}$ over $X_\Gamma$, with
the same trivialising cover as $\eta$.
Thus, if $\eta$ is numerable, so is $\eta_{\scr \Gamma}$. 
Let $(X,\pi,\varphi)$ be a \spl $\Gamma$-CW-complex over $A$.
By \thref{Th-classi-CW}, any \spl $\Gamma$-equivariant
principal bundle over $X$ is numerable.  
Hence, we get a map $\Psi\:\nbungt(X)\to[X_\Gamma,BG]$.
Also, the isotropy representation 
$\isor\:\nsbungt(X)\hfl{\approx}{}\repi$ is a bijection.
Passing to the classifying spaces gives a map
$B\:\repi\to [B\cali,BG]$.
The map $\delta\:B\cali\to X_\Gamma$ of \secref{haefborel}
gives rise to a map $\delta^*\:[X_\Gamma,BG]\to [B\cali,BG]$.

\begin{Proposition}\label{borel-bu1} 
Let $\cali$ be a \pgagrou\ for a Lie group $\Gamma$. 
Let $(X,\pi,\varphi)$ be a \spl $\Gamma$-CW-complex over $A$
with isotropy groupoid $\ciat$. 
Let $G$ be a compact Lie group.
Then, the following diagram
$$
\begin{array}{c}{\xymatrix@C-3pt@M+2pt@R-4pt{%
\nsbungt(X)\kern 3pt \ar[d]^(0.45){\isor}_(0.45){\approx} \ar[r]^{\Psi}  &
[X_\Gamma,BG] \ar[d]^(0.45){\delta^*}
\\
\repi \ar[r]^(0.50){B}  &
[B\cali,BG]
}}\end{array}
$$\nobreak
is commutative.
\end{Proposition}

\preu 
By \proref{recon-pro-CW}, we may assume that
$(X,\pi,\varphi)=(Y_{\cali},\varpi,\phi)$.
Let $\varepsilon\in\nsbungt(Y_\cali)$ and  
let $\epsilon\:\cali\to G$ be a representatative of $\isor(\varepsilon)$.
By \thref{Th-classi-CW} and its proof,
$\varepsilon$ has a representative $\eta$ of the form
$$
\Gamma\times_{\cali}(A\times G) \to
\Gamma\times_{\cali}A = Y_\cali \ ,
$$
where $\cali$ acts on $A\times G$ by
$\zeta\cdot (a,g)=(a,\epsilon(\zeta)\,g)$.
The bundle $\eta_{\scr \Gamma}$ takes the form:
$$
E\Gamma\times_{\cali}(A\times G) \to E\Gamma\times_{\cali}A = X_\Gamma \ .
$$
Let $q\:L\to B\cali$ be the induced bundle $\delta^*\eta_{\scr \Gamma}$.
To prove \proref{borel-bu1}, it is enough to construct a
$G$-equivariant map $F\:L\to EG$ making the following diagram
commutative:
$$
\xymatrix{
L \ar[r]^F\ar[d]^q&  EG\ar[d] \\
B\cali \ar[r]^{B\epsilon} & BG
}\ .
$$
Restricted to the stalk over $a$, the bundle $\delta^*\eta_{\scr \Gamma}$
is of the form
$$
E\cali_a\times_{\cali_a}(\{a\}\times G) \longrightarrow
E\cali_a\times_{\cali_a}\{a\} \ .
$$
Therefore, the required map $F$ can be defined by
\[
F(u,a,g) = E\epsilon(u)\cdot g \ . \qedhere\]
\cqfd

\proref{borel-bu1} allows us to study the map $B\:\repi\to [B\cali,BG]$,
especially when $G$ is abelian, in which case $B$ is a homomorphism of abelian groups.

\begin{Proposition}\label{gmapB} 
Let $\cali$ be a \pgagrou\ for a Lie group $\Gamma$. 
Let $(X,\pi,\varphi)$ be a \spl $\Gamma$-CW-complex over $A$
with isotropy groupoid $\ciat$. 
Let $G$ be a compact abelian Lie group.
Then, one has an isomorphism of split exact sequences of abelian groups:
\begin{equation}\label{gmapB-diag}
\begin{array}{c}{\xymatrix@C-3pt@M+2pt@R-4pt{%
0 \ar[r] &
\nsbungt(X) \ar[d]^(0.50){\isor}_(0.50){\approx}
\ar[r]  &
\nbungt(X) \ar[d]^(0.50){\delta^*\pcirc\Psi}_(0.50){\approx}
\ar[r]^(0.50){\varphi^*}  &
\bung(A) \ar[d]^(0.50){\approx} \ar@/^/[l]^{\pi^*} 
\ar[r] & 0
\\
0 \ar[r] &
\arep(\cali) \ar[r]^(0.50){B}  &
[B\cali,BG] \ar[r]^(0.50){j^*} &
[A,BG]
\ar[r] \ar@/^/[l]^{\bar\pi^*}
& 0
}}\end{array}
\ .
\end{equation}
\end{Proposition}

\preu
The top split exact sequence of abelian groups comes from \proref{Th-classi-ab-CW}
and its proof. For the bottom one, one has at least a sequence
$$\repi\xrightarrow{B} [B\cali,BG]\xrightarrow{j^*}[A,BG]$$ with
$j^*\pcirc B=0$. 
By \proref{Th-classi-ab-CW}, any principal $\Gamma$-equivariant 
$G$-bundle over $X$ is numerable. Therefore,   
the map $\delta^*\pcirc\Psi$ is defined and is a homomorphism
of abelian groups. One checks that the left-hand square
of the Diagram~\eqref{gmapB-diag} is commutative, as well as the 
right-hand square with $\varphi^*$ and $j^*$. 
The map $\delta^*$ is bijective by \proref{borel-sp}. 
As $G$ is abelian, the map $\Psi$ is a bijection by \cite[Theorem~A]{LMS}.
Thus, $\delta^*\pcirc\Psi$ is an isomorphism.   
This proves that the bottom sequence of Diagram~\eqref{gmapB-diag} is split exact.
\cqfd

\begin{Corollary}\label{mapB} 
Let $\cali$ be a \pgagrou\ for a Lie group $\Gamma$.
Let $(X,\pi,\varphi)$ be a \spl $\Gamma$-space over $A$
with isotropy groupoid $\ciat$. Let $G$ be a compact abelian Lie group.
Suppose that $H^1(A;\pi_0(G))=H^2(A;\bbz)=0$. Then the map
$B\:\repi\to [B\cali,BG]$ is a bijection.
\end{Corollary}

\preu
The abelian group $G$ is a disjoint union of tori,
so $\pi_j(BG)=\pi_{j-1}(G)=0$ for $j>2$. 
One has $\hom(\pi_1(A),\pi_1(BG))=\hom(\pi_1(A),\pi_0(G))\approx H^1(A;\pi_0(G))=0$. 
A map $f\:A\to BG$ is then null-homotopic on the $1$-skeleton and the obstruction theory
to homotop it to a constant map is with constant coefficients. Our hypotheses
implies that $H^2(A;\pi_2(BG))=0$, so one gets $[A,BG]=0$. \corref{mapB} then follows from
\proref{gmapB}.
\cqfd

\begin{ccote}\label{ktheory}
\emph{Equivariant $K$-theory}.
For vector bundles it is natural to stabilize, and then to study bundles via equivariant $K$-theory. For example, if $G =U(n)$ we consider the stablization maps
$$\sbun^{U(n)}_\Gamma(X) \to \sbun^{U(n+1)}_\Gamma(X)$$
and point out how stabilization is related to our classification results. 
\begin{Proposition}  
Let $\cali$ be a \pgagrou\ for a Lie group $\Gamma$. 
Let $(X,\pi,\varphi)$ be a \spl $\Gamma$-CW-complex over $A$
with isotropy groupoid $\cali$. Then, there
is a natural isomorphism 
$$\isor\colon K_\Gamma(X,A) \cong \krep(\cali)$$
of abelian groups induced by the isotropy representations.
\end{Proposition}
The group $\krep(\cali)$ is the Grothendieck group of the abelian monoid obtained by stabilization from the system $\{\repr^{U(n)}(\cali)\}$.
\end{ccote}
\begin{ccote}\label{eqcoh} {\it Equivariant cohomology.} \
Let $\Gamma$ be a compact Lie group and $X$ a $\Gamma$-CW-complex.
By \cite[Theorem~A]{LMS}, one has isomorphisms
\begin{equation}\label{eqcoh-eq1}
\bungs(X) \xrightarrow{\approx} [X_\Gamma,BS^1] \approx H^2_\Gamma(X) \, ,
\end{equation}
where $H^*_\Gamma(X)=H^*_\Gamma(X;\bbz)$ denotes the equivariant cohomology.
If $(X,\pi,\varphi)$ is a \spl $\Gamma$-space over a CW-complex $A$,
then the projection $\pi$ descends to a map $\bar \pi\:X_\Gamma\to A$. We denote by
$X^{(i)}$ the part of $X$ above the $i$-skeleton of $A$ and by
$r_i\:H^*_\Gamma(X)\to H^*_\Gamma(X^{(i)})$ the restriction homomorphism,
induced by the inclusion $X^{(i)}\subset X$.

\begin{Proposition}\label{P-eqcoh1} 
Let $\Gamma$ be topological group and $A$ be a CW-complex. Let $\cali$ be a \pgagrou\
and $(X,\pi,\varphi)$ be a \spl $\Gamma$-space with isotropy groupoid $\cali$.
Then
\renewcommand{\labelenumi}{(\alph{enumi})}
\begin{enumerate}
\item The sequence $0\to H^2(A)\xrightarrow{\bar\pi^*} H^2_\Gamma(X) \xrightarrow{r_0}
H^2_\Gamma(X^{(0)})$ is exact.
\item The two restriction homomorphisms
$r_0 \:H^2_\Gamma(X)\to H^2_\Gamma(X^{(0)})$ and \\ $r_{10} \:H^2_\Gamma(X^{(1)})\to H^2_\Gamma(X^{(0)})$
have the same image.
\end{enumerate}
\end{Proposition}

\begin{proof}
The map $\bar\pi$ admits a section $\bar\varphi\:A\to X_\Gamma$ coming from $\varphi$. Hence,
$\bar\pi^*\:H^*(A)\to H^*_\Gamma(X)$ is injective.

One has $H^2_\Gamma(X^{(0)})\approx \sbuns{X^{(0)}}\approx\areps(\cali^{(0)})$.
The composed homomorphism
$\buns(A)\approx H^2(A)\xrightarrow{\pi^*} H^2_\Gamma(X) \xrightarrow{r_0}
H^2_\Gamma(X^{(0)})\approx\areps(\cali^{(0)})$ sends an $S^1$-bundle $\xi$ over $A$
to the isotropy representation of $\pi^*\xi$, which is trivial.
Thus, $r_0\pcirc\bar\pi^*=0$. Using \proref{Th-classi-ab}, one has an isomorphism
$H^2_\Gamma(X)\approx \areps(\cali)\times H^2(A)$. The remainder of (a) and (b)
follow from \proref{res01skel}.
\end{proof}

We now specialise to $\Gamma$ being a torus $\bbt$, with Lie algebra $\algl$, and
use the definitions and notations of \secref{CBOBSS}.
If $\cali$ is $0$-toric, we have from Equation~\eqref{CBOBSS-eq1}, that
$$
H^2_\bbt(X^{(0)})\approx\areps(\cali^{(0)})\approx \prod_{v\in\Lambda_0} \algl^*
$$
by using \proref{P-eqcoh1} and its proof, together with \proref{imageS1}.

\begin{Proposition}\label{P-eqcph-imageS1}
Let $\cali$ be a $1$-toric \ctgrou\ and let
$(X,\pi,\varphi)$ be a \spl $\bbt$-space with isotropy groupoid $\cali$.
The image of $r\:H^2_\bbt(X)\to \prod_{v\in\Lambda_0} \algl^*$
is the set of elements $(a_v)_{v\in\Lambda_0}$
satisfying the GKM-condition.
\end{Proposition}

\begin{Remark}\label{eqcoh-rem}
Let $X$ be as in \proref{P-eqcph-imageS1}. Suppose that
$X$ is {\it equivariantly formal}, i.e.
the homomorphism $H^*_\bbt(X)\to H^*(X)$ induced by the inclusion $X\subset X_\bbt$
is surjective. In this case, the homomorphim $H^2_\Gamma(X) \xrightarrow{r_0}
H^2_\Gamma(X^{(0)})$ is injective and Part (b) of \proref{P-eqcoh1} holds,
see \cite[Theorem~1]{FP}. The injectivity of $r_0$ is considered as a
``localisation theorem'', see e.g.~\cite[Theorem~6.3]{GKM},
and Part~(b) of \proref{P-eqcoh1} is referred to as the
``Chang-Skjelbred principle'' (it historically occurred in \cite[Lemma~2.3]{CS}
for rational coefficients).
But, by Part (a) of \proref{P-eqcoh1}, $X$ is equivariantly
formal only if $H^2(A)=0$, so our context is different.
For complex coefficients,
\proref{P-eqcph-imageS1} was proven in \cite[Theorem~7.2]{GKM}. There
$X$ need not to be \spl\hskip-3pt, but again must be equivariantly formal.
\end{Remark}
\end{ccote}

\providecommand{\bysame}{\leavevmode\hbox to3em{\hrulefill}\thinspace}
\providecommand{\MR}{\relax\ifhmode\unskip\space\fi MR }
\providecommand{\MRhref}[2]{%
  \href{http://www.ams.org/mathscinet-getitem?mr=#1}{#2}
}
\providecommand{\href}[2]{#2}


\begin{thebibliography}{10}

\bibitem{Ad}
J.~F. Adams, \emph{Lectures on {L}ie groups}, W. A. Benjamin, Inc., New
  York-Amsterdam, 1969.

\bibitem{Ada}
S.~Adams, \emph{Reduction of cocycles with hyperbolic targets}, Ergodic Theory
  Dynam. Systems \textbf{16} (1996), 1111--1145.

\bibitem{Au}
M.~Audin, \emph{The topology of torus actions on symplectic manifolds},
  Progress in Mathematics, vol.~93, Birkh\"auser Verlag, Basel, 1991,
  Translated from the French by the author.

\bibitem{BH}
M.~R. Bridson and A.~Haefliger, \emph{Metric spaces of non-positive curvature},
  Grundlehren der Mathematischen Wissenschaften [Fundamental Principles of
  Mathematical Sciences], vol. 319, Springer-Verlag, Berlin, 1999.

\bibitem{CS}
T.~Chang and T.~Skjelbred, \emph{The topological {S}chur lemma and related
  results}, Ann. of Math. (2) \textbf{100} (1974), 307--321.

\bibitem{DJ}
M.~W. Davis and T.~Januszkiewicz, \emph{Convex polytopes, {C}oxeter orbifolds
  and torus actions}, Duke Math. J. \textbf{62} (1991), 417--451.

\bibitem{Du}
J.~Dugundji, \emph{Topology}, Allyn and Bacon Inc., Boston, Mass., 1978,
  Reprinting of the 1966 original, Allyn and Bacon Series in Advanced
  Mathematics.

\bibitem{DK}
J.~J. Duistermaat and J.~A.~C. Kolk, \emph{Lie groups}, Universitext,
  Springer-Verlag, Berlin, 2000.

\bibitem{FP}
M.~Franz and V.~Puppe, \emph{{Exact cohomology sequences with integral
  coefficients for torus actions}}, arXiv:math.AT/0505607.

\bibitem{GKM}
M.~Goresky, R.~Kottwitz, and R.~MacPherson, \emph{Equivariant cohomology,
  {K}oszul duality, and the localization theorem}, Invent. Math. \textbf{131}
  (1998), 25--83.

\bibitem{Gu}
V.~Guillemin, \emph{Moment maps and combinatorial invariants of {H}amiltonian
  {$T\sp n$}-spaces}, Progress in Mathematics, vol. 122, Birkh\"auser Boston
  Inc., Boston, MA, 1994.

\bibitem{Hae1}
A.~Haefliger, \emph{Homotopy and integrability}, Manifolds--Amsterdam 1970
  (Proc. Nuffic Summer School), Lecture Notes in Mathematics, Vol. 197,
  Springer, Berlin, 1971, pp.~133--163.

\bibitem{hhausmann2}
I.~Hambleton and J.-C. Hausmann, \emph{Equivariant principal bundles over
  spheres and cohomogeneity one manifolds}, Proc. London Math. Soc. (3)
  \textbf{86} (2003), 250--272.

\bibitem{KT}
Y.~Karshon and S.~Tolman, \emph{The moment map and line bundles over
  presymplectic toric manifolds}, J. Differential Geom. \textbf{38} (1993),
  465--484.

\bibitem{L2}
R.~K. Lashof, \emph{Equivariant bundles}, Illinois J. Math. \textbf{26} (1982),
  257--271.

\bibitem{L3}
\bysame, \emph{Equivariant bundles over a single orbit type}, Illinois J. Math.
  \textbf{28} (1984), 34--42.

\bibitem{LM}
R.~K. Lashof and J.~P. May, \emph{Generalized equivariant bundles}, Bull. Soc.
  Math. Belg. S\'er. A \textbf{38} (1986), 265--271 (1987).

\bibitem{LMS}
R.~K. Lashof, J.~P. May, and G.~B. Segal, \emph{Equivariant bundles with
  abelian structural group}, Proceedings of the Northwestern Homotopy Theory
  Conference (Evanston, Ill., 1982) (Providence, R.I.), Amer. Math. Soc., 1983,
  pp.~167--176.

\bibitem{Lu}
W.~L{\"u}ck, \emph{Transformation groups and algebraic {$K$}-theory}, Lecture
  Notes in Mathematics, vol. 1408, Springer-Verlag, Berlin, 1989, Mathematica
  Gottingensis.

\bibitem{Lu-O}
W.~L{\"u}ck and B.~Oliver, \emph{The completion theorem in {$K$}-theory for
  proper actions of a discrete group}, Topology \textbf{40} (2001), 585--616.

\bibitem{LW}
A.~T.~Lundell and S.~Weingram,
\emph{The topology of CW-complexes}.
Van Nostrand (1969).

\bibitem{Mo}
P.S. Mostert.
\emph{Local cross sections in locally compact groups},
Proc. Amer. Math. Soc. \textbf{4} (1953), 645--649.

\bibitem{Mo2}
P.S. Mostert.
\emph{Sections in principal fiber space},
Duke Journ. of Math. \textbf{23} (1956), 57--71.


\bibitem{St}
N. Steenrod.
 \emph{The {T}opology of {F}ibre {B}undles}.
Princeton University Press, 1951.


\bibitem{S}
T.~E. Stewart, \emph{Lifting group actions in fibre bundles}, Ann. of Math. (2)
  \textbf{74} (1961), 192--198.

\bibitem{D1}
T.~tom Dieck, \emph{Faserb\"undel mit {G}ruppenoperation}, Arch. Math. (Basel)
  \textbf{20} (1969), 136--143.

\bibitem{D2}
\bysame, \emph{Transformation groups}, de Gruyter Studies in Mathematics,
  vol.~8, Walter de Gruyter \& Co., Berlin, 1987.

\end{thebibliography}
\end{document}